\documentclass[11pt,psamsfonts]{amsart}
\usepackage{amssymb,amscd,amsmath,amsthm,latexsym}
\makeatletter \@addtoreset{equation}{section}

 \theoremstyle{plain} %% This is the default
\newtheorem{theorem}{Theorem}[section]

\newtheorem{lemma}[theorem]{Lemma}
\newtheorem{proposition}[theorem]{Proposition}

\theoremstyle{remark}
\newtheorem{remark}[theorem]{Remark}
\newtheorem{example}[theorem]{Example}
\theoremstyle{definition}
\newtheorem{definition}[theorem]{Definition}

\newcommand{\thmref}[1]{Theorem~\ref{#1}}

\newcommand{\proref}[1]{Proposition~\ref{#1}}
\newcommand{\lemref}[1]{Lemma~\ref{#1}}

\def\1{{\bf 1}}

\def\2{\frac{1}{2}}

\def\C{{\mathbb C}}
\def\NN{{\mathbb N}}
\def\R{{\mathbb R}}

\def\Z{{\mathbb Z}}

\def\F{{\mathcal F}}
\def\L{{\mathcal L}}
\def\M{{\mathcal M}}
\def\N{{\mathcal N}}

\def\O{{{\mathcal O}_X}}

\def\To{{{\mathcal T}_X}}
\def\U{{\mathcal U}}
\def\sa{_\frac{i\beta}{2}}
\def\sb{_{-\frac{i\beta}{2}}}
\def\sc{_\frac{i}{2}}
\def\sd{_{-\frac{i}{2}}}

\def\eps{\varepsilon}

\def\Ad{{\rm Ad}\,}

\def\dim{{\rm dim}\,}
\def\End{{\rm End}}
\def\Hom{{\rm Hom}}

\def\Ker{{\rm Ker}}

\def\Mat{{\rm Mat}}

\def\Sp{{\rm Sp}}
\def\span{{\rm span}}

\def\supp{{\rm supp}\,}
\def\Tr{{\rm Tr}_\tau}
\def\Trr{{\rm Tr}_{\tau_0}}
\def\Ttr{{\rm Tr}}
\def\<{\langle}
\def\>{\rangle}

\begin{document}

\title[KMS states of quasi-free dynamics ]{KMS states of quasi-free dynamics on Pimsner algebras}

\author[Laca]{Marcelo Laca*}
\address{Department of Mathematics and Statistics,
University of Victoria,
Victoria,  V8W 3P4 Canada.}
\email{laca@math.uvic.ca}
\thanks{* Supported by the National Science and Engineering Research
Council of Canada.}

\author[Neshveyev]{Sergey Neshveyev**}
\address{Mathematics Institute,
University of Oslo,
PB 1053 Blindern,
Oslo 0316,
Norway.}
\email{neshveyev@hotmail.com}
\thanks{** Supported by
the Royal Society/NATO postdoctoral fellowship, the Centre for Advanced
Study in Oslo and the Norwegian Research Council.}

\date{April 28, 2003}

\begin{abstract}  
A continuous one-parameter group of unitary isometries of a right-Hilbert C*-bimodule induces 
a quasi-free dynamics on the Cuntz-Pimsner C*-algebra of the bimodule and on its Toeplitz extension. 
The restriction of such a dynamics to the algebra of coefficients of the bimodule is trivial, 
and the corresponding KMS states of the Toeplitz-Cuntz-Pimsner and Cuntz-Pimsner C*-algebras 
are characterized in terms of traces on
the algebra of coefficients. This generalizes and sheds light onto
various earlier results about KMS states of the gauge actions on Cuntz
algebras,  Cuntz-Krieger algebras, and crossed products by
endomorphisms. We also obtain a more general characterization, in terms of KMS weights,
for the case in which the inducing isometries are not unitary, and accordingly, 
the restriction of the quasi-free dynamics to the algebra of coefficients is nontrivial.

\end{abstract}
\maketitle

\section*{Introduction}
Soon after the introduction of the Cuntz algebras in \cite{C} it
was noticed that the gauge action on $\mathcal O_n$ had the unique
equilibrium inverse temperature $\beta = \log n$,
\cite{OP,Ev,BEK}. Along the same lines, the gauge action on the
Cuntz-Krieger algebra $\mathcal O_A$ was also shown
 to have a unique KMS state, at inverse temperature equal to the
logarithm of the spectral radius of the irreducible matrix $A$,
\cite{EFW}. Here we are interested in the KMS states of the
C*-algebras $\To$ and $\O$ associated by Pimsner to a right
Hilbert bimodule $X$, \cite{Pim}. Having a rich, yet tractable
structure, they provide a convenient framework in which to study
the interesting phenomena that characterize the
examples mentioned above and many others; see e.g. \cite{PWY} and \cite{EL2}.
Specifically, we start with a C$^*$-algebra $A$ and a right
Hilbert $A$ bimodule $X$ in which the left action is
non-degenerate. Given a continuous one-parameter group of
isometries on $X$, we induce quasi-free dynamics on $\To$ and on
$\O$ via their universal properties, cf \cite{Z}. We then proceed
to study the equilibrium states of these quasi-free dynamics
associated to groups of isometries in terms of their restrictions
to the coefficient algebra $A$. Our approach underlines the role
of the Toeplitz algebra $\To$ in its own right and not as a mere
preliminary step from which to obtain $\O$ as a quotient. The key
point, inspired in Evans construction \cite{Ev}, is that $\To$
acts naturally on the full Fock module over $X$, and the
quasi-free action is implemented there by the Fock quantization of
the given group of isometries on $X$. We use this as a guidance in
writing KMS states as quasi--free states, but do not rely on it
directly in our arguments. Since $\O$ does not act in general on
the Fock module, this type of spatial (modular) implementation is
lost when one looks at quasi-free dynamics on  $\O$ alone.
However, it is easy to characterize the KMS states on $\O$ as
those on $\To$ that factor through the quotient.

A brief summary of the contents follows. In Section 1 we collect
some necessary results about inducing traces from the coefficient
algebra $A$ to the algebra of adjointable operators on $X$. In
Section 2, after introducing the Pimsner algebra of a bimodule and
the quasi-free dynamics associated to a one-parameter group of
isometries of the bimodule, we study the special case of dynamics
that fix the elements of $A$. Under a positivity assumption on the
infinitesimal generator $D$ of the given group of isometries, we
show that the KMS$_\beta$ states of the quasi-free dynamics are
induced from the traces on $A$ that satisfy a certain inequality.
This inequality is formulated in terms of a transfer operator between
traces (or KMS weights) on $A$
 and, essentially, ensures that the Fock quantization of the
contraction $e^{-\beta D}$ is an appropriate density operator. In
the special case of the gauge action on the Toeplitz Cuntz
algebras, the coefficient algebra is $\mathbb C$, and its unique
trace satisfies the inequality if and only if $\beta \geq \log n$;
the resulting KMS$_\beta$ state of $\mathcal T_n$ factors through
$\mathcal O_n$ only for $\beta = \log n$. At the end of the
section we show how to derive, in a unified way, several other
examples of KMS states of C*-dynamical systems previously studied
under different guises. In Section 3 we extend our results to the
more general situation in which the dynamics on the coefficient
algebra $A$ is allowed to be nontrivial. Most of the section is
devoted to inducing KMS states on $A$. There is no great
simplification at this point in restricting ourselves to states,
so in fact we consider KMS weights on $A$. For von
Neumann algebras this problem was studied  by Combes and Zettl in
\cite[Section 3]{CZ}, who deduced the existence of induced weights
from the well-known cocycle theorem of Connes. Besides giving a
slightly different contruction which works equally well for
C$^*$-algebras, we provide also a direct proof and indicate how
Connes' result can be derived from our results on induced weights.
Once the induction procedure for weights has been settled, the
characterization of KMS states in terms of their restrictions to
$A$ is entirely analogous to that of KMS states in Section 2.

This work was initiated during a short visit of M.L. to Cardiff
University and continued through visits, of M.L. to the Center for
Advanced Studies in Oslo, and of S.N. to the University of
Victoria. The authors would like to acknowledge the support and
the hospitality provided by those institutions.

\section{Preliminary results on induced traces.}

If $A$ is an algebra and $X$ is a right projective $A$-module of
finite type, then the algebra $\End_A(X)$ is isomorphic to
$X\otimes_A\Hom_A(X,A)$. Hence there exists a unique linear map
$\Ttr\colon\End_A(X)\to A/[A,A]$ such that $\Ttr(x\otimes
f)=f(x)\,{\rm mod}\,[A,A]$. The composition of any tracial linear
functional~$\tau$ on $A$ (one for which $\tau(ab)=\tau(ba)$) with
$\Ttr$ yields an (induced) tracial linear functional $\Tr$ on
$\End_A(X)$. Clearly, this construction can be applied to any
unital C$^*$-algebra $A$ and finite Hilbert $A$-module $X$. Our
first aim in this section is to define $\Tr$ for arbitrary
C$^*$-algebras and Hilbert modules, and to derive some of its
basic properties.

Suppose $X$ is a right Hilbert module over a C$^*$-algebra $A$ and
let $K(X)$ be the C$^*$-algebra of generalized compact operators
on $X$, generated by the operators $\theta_{\xi,\zeta}$, given by
$\theta_{\xi,\zeta}\eta=\xi\<\zeta,\eta\>$, with $\xi,\zeta, \eta
\in X$. Let $B(X)$ be its multiplier algebra, that, is, the
C*-algebra of adjointable operators on $X$. Recall that a bounded
net $\{S_k\}_k$ in $B(X)$ converges to $S\in B(X)$ strictly if and
only if $S_k\xi\to S\xi$ and $S^*_k\xi\to S^*\xi$ for every
$\xi\in X$. We shall need the following result about inducing
traces from $A$ through $X$. The existence and some of the
properties of the induced trace $\Tr$ can be found in
\cite[Section 2]{CZ}, and \cite[Lemma 4.6]{CPtraces}, where they
are derived from previous results of Pedersen \cite{P} about
extending traces from hereditary subalgebras. We state the
relevant properties in a way that is convenient for our purposes,
and we include a self-contained, direct proof for completeness.
 \begin{theorem} \label{1.1}
Let $\tau$ be a finite trace on $A$. For $T\in B(X)$, $T\ge0$, set
$$
\Tr(T)=\sup_I\sum_{\xi\in I}\tau(\<\xi,T\xi\>),
$$
where the supremum is taken over all finite subsets $I$ of
$X$ such that
$\sum_{\xi\in I}\theta_{\xi,\xi}\le1$. Define
\[
\M^+_\tau=\{T\ge0\,|\,\Tr(T)<\infty\},\quad
\N_\tau=\{T\,|\,T^*T\in\M^+_\tau\}, \quad
\M_\tau=\span\,\M^+_\tau=\N^*_\tau\N_\tau.
\]\
Then
\begin{enumerate}
\item[\rm(i)] $\Tr$ is strictly lower semicontinuous,
moreover, if $\liminf_k\tau(\<\xi,T_k\xi\>)\ge\tau(\<\xi,T\xi\>)$
for every $\xi\in X$, then $\liminf_k\Tr(T_k)\ge\Tr(T)$;

\smallskip
\item[\rm(ii)] if $\{e_k=\sum_{\xi\in I_k}\theta_{\xi,\xi}\}_k$ is
a net such that $e_k\le1$ and
$\tau(\<\eta,e_k\eta\>)\to\tau(\<\eta,\eta\>)$ for every $\eta\in
X$ (e.g. if $\{e_k\}_k$ is an approximate unit in $K(X)$), then
$\Tr(T)=\lim_k\sum_{\xi\in I_k}\tau(\<\xi,T\xi\>)$ for $T\in
B(X)_+$; in particular, $\Tr$ can be extended to a positive linear
functional on $\M_\tau$;

\smallskip
\item[\rm(iii)] for every pair $\xi,\eta\in X$, we have that
$\theta_{\xi,\eta}\in{\M}_\tau$ and
$\Tr(\theta_{\xi,\eta})=\tau(\<\eta,\xi\>)$;

\smallskip
\item[\rm(iv)] $\Tr$ is a semifinite trace; thus $\N_\tau$ and
$\M_\tau$ are two-sided ideals in $B(X)$, $\M_\tau$ is essential,
and if $S,T\in\N_\tau$, or if $S\in B(X)$ and $T\in\M_\tau$, then
$\Tr(ST)=\Tr(TS)$.
\end{enumerate}
\end{theorem}

\begin{proof}
The proof of (i) is trivial. To prove (ii) we shall first
prove that $\Tr(\theta_{\xi,\xi})=\tau(\<\xi,\xi\>)$. Suppose
$S=\sum_{\eta\in I}\theta_{\eta,\eta}\le1$. Then
$$
\sum_{\eta\in I}\tau(\<\eta,\theta_{\xi,\xi}\eta\>)
=\sum_{\eta\in I}\tau(\<\eta,\xi\>\<\xi,\eta\>)
=\sum_{\eta\in I}\tau(\<\xi,\eta\>\<\eta,\xi\>)=\tau(\<\xi,S\xi\>),
$$
and the equality $\Tr(\theta_{\xi,\xi})=\tau(\<\xi,\xi\>)$
follows. The same proof shows that $\Tr(\sum_{\xi\in
I}\theta_{\xi,\xi})=\sum_{\xi\in I}\tau(\<\xi,\xi\>)$ for any
finite set $I$. Hence if $\{e_k\}_k$ is as in the formulation of
(ii) and $T\in B(X)_+$, then
$$
\Tr(T^\2e_kT^\2)=\Tr(\sum_{\xi\in I_k}
\theta_{T^\2\xi,T^\2\xi}) =\sum_{\xi\in
I_k}\tau(\<\xi,T\xi\>).
$$
Since $\lim_k\Tr(T^\2e_kT^\2)=\Tr(T)$ by property~(i), (ii)
follows. Part (iii) has already been proved for $\xi=\eta$, the
general case follows by polarization, and implies that $\M_\tau$
is essential. Since $\theta_{u\xi,u\xi}=u\theta_{\xi,\xi}u^*$, it
is obvious that $\Tr(uTu^*)=\Tr(T)$ for any unitary $u\in B(X)$.
Thus $\Tr$ is a trace.
\end{proof}

Suppose now $Y$ is a right
Hilbert $A$-bimodule, that is, $Y$ is a right Hilbert $A$-module
together with a left action of $A$ given by a $*$-homomorphism of
$A$ into $B(Y)$. Denote by $B_A(Y)$ the subalgebra of $B(Y)$
consisting of $A$-bimodule maps. Let $X$ be another right Hilbert
$A$-module. Then the tensor product $X\otimes_AY$ is a right
Hilbert $A$-module and for any $S\in B(X)$ and $T\in B_A(Y)$ there
is an operator $S\otimes T\in B(X\otimes_AY)$.
We shall need the following property about
the induction in stages of a trace $\tau$ on $A$ through a tensor
product of modules. We indicate explicitly in the proposition the
module used to induce the trace in each case, but we drop this
notation later for simplicity, and rely on the context to indicate
the relevant module.

\begin{proposition} \label{1.2}
Let $S\in B(X)_+$, and $T\in B_A(Y)_+$. Suppose $\tau$ is a finite
trace on $A$ such that $\Tr^Y(T)<\infty$, and define a new finite
trace $\tau_T$ on $A$ by letting $\tau_T(a)=\Tr^Y(aT)$. Then
$\Tr^{X\otimes Y}(S\otimes T)=\Ttr_{\tau_T}^X(S)$.
 \end{proposition}

\begin{proof} We begin by constructing an
approximate unit for $K(X\otimes_AY)$ from approximate units for
$K(X)$ and $K(Y)$. For finite subsets $I \subset X$ and $J \subset
Y$, define the operators $e_I = \sum_{\xi\in I}\theta_{\xi,\xi}$,
$e_J= \sum_{\zeta\in J}\theta_{\zeta,\zeta}$ and
$e_{I,J}=\sum_{\xi\in I,\zeta\in
J}\theta_{\xi\otimes\zeta,\xi\otimes\zeta}$. We claim that if
$e_I\le1$ and $e_J\le1$, then $e_{I,J}\le1$. To verify this we
consider a vector $\eta=\sum_k\mu_k\otimes\nu_k$ in $X\otimes_AY$;
and observe that
$$
\<e_{I,J}\eta,\eta\>=\sum_{\xi\in I}\<e_J\delta_\xi,\delta_\xi\>,
$$
where $\delta_\xi=\sum_k\<\xi,\mu_k\>\nu_k$. If $e_I\le1$, then
$(\<e_I\mu_k,\mu_l\>)_{k,l}\le(\<\mu_k,\mu_l\>)_{k,l}$ in the
algebra $\Mat_{|I|}(A)$ of $|I|$ by $|I|$ matrices over $A$, so if,
in addition, $e_J\le1$ we get
$$
\sum_{\xi\in I}\<e_J\delta_\xi,\delta_\xi\>
\le\sum_{\xi\in I}\<\delta_\xi,\delta_\xi\>
=\sum_{k,l}\<\nu_k,\<e_I\mu_k,\mu_l\>\nu_l\>
\le\sum_{k,l}\<\nu_k,\<\mu_k,\mu_l\>\nu_l\>
=\<\eta,\eta\>,
$$
which proves the claim. Note also that
$$
e_{I,J}(\mu\otimes\nu)=e_I\mu\otimes\nu
+\sum_{\xi\in I}\xi\otimes(e_J-1)\<\xi,\mu\>\nu.
$$
It follows that there exists an approximate unit in
$K(X\otimes_AY)$ consisting of elements of the form~$e_{I,J}$, with
$e_I\le1$ and $e_J\le1$, so, by \thmref{1.1}(ii),
$$
\Tr^{X\otimes Y}(S\otimes T)=\sup_{I,J}\sum_{\xi\in I,\zeta\in J}
\tau(\<\zeta,\<\xi,S\xi\>T\zeta\>).
$$
Since for fixed $I$,
$$
\sup_J\sum_{\xi\in I,\zeta\in J}\tau (\<\zeta,\<\xi,S\xi\>T\zeta\>)
=\sum_{\xi\in I}\Tr^Y(\<\xi,S\xi\>T)=\sum_{\xi\in
I}\tau_T(\<\xi,S\xi\>),
$$
recalling that $\displaystyle\sup_I\sum_{\xi\in
I}\tau_T(\<\xi,S\xi\>)=\Ttr_{\tau_T}^X(S)$, we get $\Tr^{X\otimes
Y}(S\otimes T)=\Ttr_{\tau_T}^X(S)$.
\end{proof}

We have shown that any finite trace on $A$ can be induced to a
unique, strictly densely defined, strictly lower semicontinuous
trace $\Tr$ on $B(X)$ such that
$\Tr(\theta_{\xi,\eta})=\tau(\<\eta,\xi\>)$. Clearly one should
not expect all strictly densely defined, strictly lower
semicontinuous traces on $B(X)$ to be induced from finite traces
on $A$. In Section~\ref{kmsweights} below we generalize this
extension procedure so that it applies to KMS weights of
quasi-free dynamics. By setting the dynamics to be trivial, we
then obtain, as a corollary, a bijective correspondence between
densely defined lower semicontinuous traces on
$\overline{\<X,X\>}$ and strictly densely defined, strictly lower
semicontinuous traces on $B(X)$.

\section{KMS states on Pimsner algebras.}

Next we consider a Hilbert $A$-bimodule $X$, with the purpose of
studying KMS-states on the Toeplitz-Pimsner and
Cuntz-Krieger-Pimsner algebras associated in \cite{Pim} to such a
bimodule. The only extra assumption that we make on the bimodule
is the non-degeneracy of the left action, i.e. $AX=X$. In
particular, we do not assume that $X$ is full or that $A$ is
unital. We denote by $i_X\colon A\to B(X)$ the homomorphism
defining the left action of $A$ on $X$.

Let ${\F}(X)=A\oplus X\oplus(X\otimes_A X)\oplus\ldots$ be the
Fock Hilbert bimodule of $X$. The Toeplitz-Pimsner algebra $\To$
of $X$ is, by definition, the C*-algebra of operators on ${\F}(X)$
generated by the left multiplication operators $i_{\F(X)}(a)$ for
$a\in A$ and the left creation operators $T_\xi$ for $\xi\in X$,
which are given by $T_\xi(\xi_1\otimes\ldots\otimes\xi_n)
=\xi\otimes\xi_1\otimes\ldots\otimes\xi_n$. It is shown in
\cite{Pim,FR} that $\To$ is the universal C$^*$-algebra generated
by elements $\pi(a)$ with $a\in A$ and $T_\xi$ with $\xi\in X$,
such that $\pi\colon A\to\To$ is a $*$-homomorphism,
$X\ni\xi\mapsto T_\xi$ is an $A$-bilinear map, (that is, $T_{\xi
a}=T_\xi\pi(a)$ and $T_{a\xi}=\pi(a)T_\xi$), and $T^*_\xi
T_\zeta=\pi(\<\xi,\zeta\>)$. More precisely, the Fock realization
of these relations, given by the left action of $A$ and the left
creation operators on $\F(X)$, determines an isomorphism of the
universal C$^*$-algebra of the relations onto $\To$.

Let $j_X\colon K(X)\to\To$ be the injective homomorphism given by
$j_X(\theta_{\xi,\eta})=T_\xi T^*_\eta$. Let also $I_X$ be the
ideal in $A$ consisting of elements $a\in A$ such that $i_X(a)\in
K(X)$. The Cuntz-Krieger-Pimsner algebra $\O$ is, by
definition, the quotient of $\To$ by the ideal generated by elements
of the form $\pi(a)-(j_X\circ i_X)(a)$, for $a\in I_X$.
We shall usually omit $\pi$ and $i_X$ in the computations below.

Let $\R\ni t\to\sigma_t$ be a one-parameter automorphism group of
$A$ and let $\R\ni t\mapsto U_t$ be a one-parameter group of
isometries on $X$ such that $U_ta\xi=\sigma_t(a)U_t\xi$ and
$\<U_t\xi,U_t\zeta\>=\sigma_t(\<\xi,\zeta\>)$; as usual, both
$\sigma$ and $U$ are assumed to be strongly continuous. By the
universal property of the Toeplitz-Pimsner algebra there exists,
for each $t\in \R$, a unique automorphism $\gamma_t$ of $\To$ such
that $\gamma_t(a)=\sigma_t(a)$ and $\gamma_t(T_\xi)=T_{U_t\xi}$.
The resulting one-parameter group $t \mapsto \gamma_t$ is strongly
continuous and is called the {\em quasi-free dynamics} associated
to the module dynamics $U$. Since $A$ and $j_X(K(X))$ are
invariant under $\gamma_t$, so is $I_X$, and thus there is a
quasi-free dynamics at the level of $\O$, too. When we view $\To$
as acting on the Fock bimodule, the automorphisms $\gamma_t$ are
implemented by the 'Fock quantization' of the isometries $U_t$;
specifically $\gamma_t = \Ad \Gamma(U_t)$, where
$\Gamma(U_t)=1\oplus U_t\oplus(U_t\otimes U_t)\oplus\ldots$. The
group of gauge transformations is a particular case of this,
corresponding to the trivial dynamics on $A$ and the one-parameter
scalar unitary group $\{\xi\mapsto e^{it}\xi\}_{t\in\R}$ on $X$.

Given a quasi-free dynamics on $\To$, we are interested in the
relation between the $(\sigma,\beta)$-KMS states on $A$ and the
$(\gamma,\beta)$-KMS states on $\To$. For simplicity we shall
first consider this question under the assumption that the
dynamics $\sigma$ on $A$ is trivial. This covers most examples in
the literature, and has the advantage of being tractable using the
elementary properties of induced traces from the preceding
section. The case of nontrivial $\sigma$ requires a generalization
of these properties to induced KMS-weights, a task that we take up
in the next section. Accordingly, we now restrict our attention to
one-parameter groups of isometries $U_t$ of $X$ such that $U_t
a\xi=aU_t\xi$ and $\<U_t\xi,U_t\zeta\>=\<\xi,\zeta\>$, in other
words, we assume that $U$ is a one-parameter {\em unitary} group
in $B_A(X)$.

\begin{theorem} \label{old3.1}\label{mainKMSb}
Let $\R\in t\to U_t$ be a one-parameter unitary group in $B_A(X)$
satisfying the following 'positive energy' condition: the vectors
$\xi\in X$ such that $\Sp_U(\xi)\subset(0,+\infty)$ form a dense
subspace of $X$, where $\Sp_U(\xi)$ is the Arveson spectrum of
$\xi$ with respect to $U$. Let $\gamma$ be the corresponding
dynamics on $\To$, given by $\gamma_t(T_\xi)=T_{U_t\xi}$ and
$\gamma_t(a)=a$, and suppose $\beta \in(0,\infty)$. If $\phi$ is a
$(\gamma,\beta)$-KMS state on $\To$, then $\tau=\phi|_A$ is a
tracial state on $A$ and
\begin{equation}\label{ineqtrace}
\Tr(ae^{-\beta D})\le\tau(a)\quad \text{for }\ a\in A_+,
\end{equation}
where $D$ is the generator of $U$ (so that $U_t=e^{itD}$).
Conversely, if $\tau$ is a tracial state on $A$ such that
\eqref{ineqtrace} is satisfied, then there exists a unique
$(\gamma,\beta)$-KMS state $\phi$ on $\To$ with $\phi|_A=\tau$.
The state $\phi$ is determined by $\tau$ through
\begin{equation}\label{2PtFunction}
\phi(T_{\xi_1} \cdots T_{\xi_m} T_{\eta_n}^* \cdots T_{\eta_1}^*)
= \left\{\begin{array}{ll} \tau ( \< \eta_1\otimes \cdots \otimes
\eta_n \, , \, e^{-\beta D}\xi_1\otimes \cdots \otimes e^{-\beta
D}\xi_n\>)
& \text{ if } m=n\\
0 & \text{ otherwise}.
\end{array}
\right.
\end{equation}
\end{theorem}

\begin{proof}
Note that our positive energy condition is equivalent to the
existence of an increasing sequence of $U$-invariant submodules
$Y_n$ of $X$ having dense union and such that $D|_{Y_n}\ge c_n 1$
with $c_n >0$. It follows from this that $e^{-\beta D}$ is a
selfadjoint contraction in $B_A(X)$ for each $\beta > 0$.

Assume first that $\phi$ is a $(\gamma,\beta)$-KMS state. Since
the left action of $A$ on $X$ is non-degenerate, the homomorphism
$\pi\colon A\to\To$ is non-degenerate. Hence $\tau=\phi|_A$ is a
state. Since $\gamma$ is trivial on $A$, $\tau$ is a trace.

For any $\xi\in X$ we have $\gamma\sa(T_\xi)=T_{e^{-\frac{\beta
D}{2}}\xi}$, so by the KMS condition we get
$$
\phi(T_\xi T^*_\eta)=\phi(\gamma\sa(T_\eta)^*\gamma\sa(T_\xi))
=\tau(\<e^{-\frac{\beta D}{2}}\eta,e^{-\frac{\beta D}{2}}\xi\>)
=\Tr(\theta_{\xi,\eta}e^{-\beta D}).
$$
Thus, for $a\in A_+$ we have
$$
\Tr(a^\2\theta_{\xi,\xi}a^\2e^{-\beta D})
=\Tr(\theta_{a^\2\xi,a^\2\xi}e^{-\beta D})
=\phi(T_{a^\2\xi}T^*_{a^\2\xi})=\phi(a^\2j_X(\theta_{\xi,\xi})a^\2).
$$
Since $\phi(a^\2j_X(\cdot)a^\2)$ is a positive linear functional
on $K(X)$ of norm\,$\le\tau(a)$, by the strict lower
semicontinuity of $\Tr$ we conclude that $\Tr(ae^{-\beta
D})\le\tau(a)$.

Let us now prove that $\phi$ is completely determined by $\tau$.
First note that there is an $A$-bilinear isometry of $X^{\otimes
n}$ to $\To$ mapping $\xi=\xi_1\otimes_A\ldots\otimes_A\xi_n$ to
$T_\xi=T_{\xi_1}\ldots T_{\xi_n}$, so that $A$ and the elements of
the form $T_\xi T^*_\zeta$, with $\xi\in X^{\otimes n}$ and
$\zeta\in X^{\otimes m}$, span a dense subspace of $\To$. The same
computation as above shows that for $\xi,\zeta\in X^{\otimes n}$
we have $\phi(T_\xi T^*_\zeta)=\tau(\<(e^{-\frac{\beta
D}{2}})^{\otimes n}\zeta,(e^{-\frac{\beta D}{2}})^{\otimes
n}\xi\>)$. So, in order to prove that $\phi$ is uniquely
determined by $\tau$ and that \eqref{2PtFunction} holds, it
suffices to show that $\phi(T_\xi T^*_\zeta)=0$ when $\xi\in
X^{\otimes n}$, $\zeta\in X^{\otimes m}$ and $n\ne m$. We may
assume that $n>m$ and $\xi=\xi_1\otimes\xi_2$ with $\xi_1\in
X^{\otimes m}$ and $\xi_2\in X^{\otimes(n-m)}$. Using the KMS
condition we get
$$
\phi(T_\xi T^*_\zeta)=\phi(T_{\xi_1}T_{\xi_2}T^*_\zeta)
=\phi(T_{\xi_2}T^*_\zeta T_{(e^{-\beta D})^{\otimes m}\xi_1})
=\phi(T_{\xi_2\<\zeta,(e^{-\beta D})^{\otimes m}\xi_1\>}).
$$
Thus, it is enough to show that $\phi(T_\xi)=0$ for every $\xi\in
X^{\otimes n}$ with $n\ge1$. We claim that the elements of the
form $\eta - (U_t)^{\otimes n}\eta$ span a dense subspace of
$X^{\otimes n}$; this will finish the proof because $\phi(T_{\eta
- (U_t)^{\otimes n}\eta})=0$ by virtue of the
$\gamma_t$-invariance of $\phi$. To prove the claim, suppose that
$\xi\in X$ is such that $\Sp_U(\xi)$ is compact and
$0\notin\Sp_U(\xi)$, and notice that such elements are dense
because of our assumption on the spectrum of $U$. Choose
$t_0\in\R$ such that the function $t \mapsto 1-e^{itt_0}$ is
non-zero on $\Sp_U(\xi)$; then $\xi=(1-U_{t_0})U_f\xi$ for every
function $f\in L^1(\R)$ such that $\hat f(t)=(1-e^{itt_0})^{-1}$
for $t$ in a neighbourhood of $\Sp_U(\xi)$. Hence $\xi$ is of the
form $\eta-U_t\eta$, proving the claim for $n = 1$. Since
$\Sp_{U^{\otimes n}}(\xi_1\otimes\ldots\otimes\xi_n)
\subset\overline{\Sp_U(\xi_1)+\ldots+\Sp_U(\xi_n)}$ by \cite{arv},
$U^{\otimes n}$ satisfies the same spectral assumption as $U$, and
the above argument also proves the claim for $n >1$. This finishes
the proof that $\phi$ is determined by $\tau$.

Denote by $F$ the operator, mapping finite traces on $A$ into
possibly infinite traces on $A$, defined by
$(F\tau)(a)=\Tr(ae^{-\beta D})$. The second part of the theorem
says that if $F\tau\le\tau$ for a tracial state~$\tau$ on $A$,
then there exists a $(\gamma,\beta)$-KMS state $\phi$ on $\To$
such that $\phi|_A=\tau$. Suppose for a moment that the tracial
state $\tau$ is of the form $\tau=\sum^\infty_{n=0}F^n\tau_0$ for
some finite trace $\tau_0$ (such a state clearly satisfies
$F\tau\le\tau$, in fact $F^n \tau \searrow 0$). We claim that in
this case the extension is given by $\Phi=\Trr(\cdot\,
\Gamma(e^{-\beta D}))$, where $\Gamma(e^{-\beta D}) = \sum_n
(e^{-\beta D})^{\otimes n}$ is the operator on $\F(X)$ obtained by
'Fock quantization' of the contraction $e^{-\beta D}$. Indeed, it
is easy to see that $\Phi$ is a positive linear functional with
the KMS property, but one must still verify that $\Phi$ is a state
extending~$\tau$. Using Proposition~\ref{1.2} we see by induction
that $\Ttr_{\tau_0}(\cdot \, (e^{-\beta D})^{\otimes
n})|_A=F^n\tau_0$, whence $\Phi|_A=\tau$. Since the left action of
$A$ on $\F(X)$ is non-degenerate, and $\Phi$ is strictly lower
semicontinuous by \thmref{1.1}, this implies that $\Phi$ is a
state. Thus $\phi=\Phi|_\To$ is the required $(\gamma,\beta)$-KMS
state extending $\tau=\sum^\infty_{n=0}F^n\tau_0$.

Suppose now that $\tau$ is an arbitrary tracial state such that
$F\tau\le\tau$. For each $\eps>0$ consider a one-parameter unitary
group $U^\eps$ defined by $U^\eps_t\xi=e^{i\eps t}U_t\xi$. Let
$\gamma^\eps$ be the associated quasi-free dynamics on $\To$. For
the corresponding operator $F_\eps$ on traces of $A$ we have
$F_\eps=e^{-\beta\eps}F$. In particular, $F_\eps\tau\le
e^{-\beta\eps}\tau$. Then we may write
$\tau=\sum^\infty_{n=0}F^n_\eps\tau_\eps$, with
$\tau_\eps=\tau-F_\eps\tau$; indeed, since $F^m_\eps\tau\le
e^{-\beta\eps m}\tau\to0$ as $m\to\infty$, we have
$\sum^m_{n=0}F^n_\eps\tau_\eps=\tau-F^{m+1}_\eps\tau \to \tau$.
Hence there exists a $(\gamma^\eps,\beta)$-KMS state $\phi_\eps$
on~$\To$ such that $\phi_\eps|_A=\tau$. As in \cite[Proposition
5.3.25]{bra-rob}, any weak* limit point of the states $\phi_\eps$
as $\eps\to0^+$ is a $(\gamma,\beta)$-KMS state $\phi$ extending
$\tau$.
\end{proof}

The situation for ground states is slightly different, since for
$\beta = \infty$ there is no tracial condition on $A$. For each
state $\omega$ of $A$ we define a {\em generalized Fock state}
$\phi_\omega$ of $\To$ by $\phi_\omega(T) = \lim_\lambda
\omega(\<e_\lambda, T e_\lambda\>)$ for $T\in \To$; where
$(e_\lambda)_{\lambda \in \Lambda}$ is an approximate unit in $A$.
Clearly $\omega = \phi_\omega|_A$ and $\phi_\omega$ is
characterized by $\phi_\omega (a) = \omega(a)$ for $a\in A $, and
$ \phi_\omega(T_\xi T_\eta^*) = 0$ for $\xi \in X^{\otimes m}$,
and $\eta \in X^{\otimes n}$ with $m$ or $n$ nonzero.

\begin{theorem}
Under the assumptions of \thmref{mainKMSb}, a state of $\To$ is a
ground state for $\gamma$, if and only if it is a generalized Fock
state.
\end{theorem}
\begin{proof}
Suppose first that $\phi$ is a generalized Fock state of $\To$. To
see that whenever $b \in \To$, $\xi \in X^{\otimes m}$ and $\eta
\in X^{\otimes n}$, the analytic function $ z \mapsto \phi(b
\gamma_z(T^{}_\xi T^*_{\eta}))$ is bounded on the upper half
plane, we write
\[|\phi(b \gamma_z(T^{}_\xi T^*_{\eta}))|
= |\phi( b T^{}_{U_z\xi} T^*_{U_{\overline{z}} \eta})| \leq
\phi(bb^*) \phi( T^{}_{U_{\overline{z}} \eta} T^*_{U_z\xi}
T^{}_{U_z\xi} T^*_{U_{\overline{z}} \eta});
\]
the right hand side vanishes for $n>0$ because $
\phi(T_{U_{\overline{z}} \eta} T_{U_{\overline{z}} \eta}^*) = 0$,
and it is bounded for $n=0$ because $\|U_z\xi \|\leq \|\xi\|$ for
$z$ in the upper half plane. Hence $\phi$ is a ground state.

Suppose next that $\phi$ is a ground state. By
\cite[Proposition 5.3.19(4)]{bra-rob},
\[
\phi( \gamma_f(T)^* \gamma_f(T)) = 0, \qquad (T \in \To)
\]
for every function $f \in L^1(\mathbb R)$ such that $\supp \hat{f}
\subset (-\infty, 0)$. By the positivity condition, the set of all
$\xi \in X$ such that $\Sp_U (\xi) $ is a compact subset of $(0,
\infty) $ is dense in $X$. For each such $\xi$ one has that
 $\Sp_\gamma(T_\xi^*)$ is a compact subset of $ (-\infty, 0)$, and there
exists a function $f$ as above such that $\gamma_f(T^*_\xi) =
T^*_\xi$. Putting $T = T^*_\xi$ one sees that  $\phi(T_\xi
T^*_\xi) = 0$ for every $\xi$ in a dense set and hence for all
$\xi \in X$. Since the one-parameter group $U^{\otimes n}$ on
$X^{\otimes n}$ satisfies the same positivity condition, this
implies that $\phi$ is a generalized Fock state.
\end{proof}

Suppose $\gamma = \Ad\Gamma(U)$ is a quasi-free dynamics
satisfying the hypothesis of \thmref{old3.1}, let $\phi$ be a
$(\gamma,\beta)$-KMS state on $\To$, and let $\tau = \phi|_A$. In
the course of the proof of \thmref{old3.1} we let
$(F\tau)(a)=\Tr(ae^{-\beta D})$ and showed that if
$\tau=\sum^\infty_{n=0}F^n\tau_0$ for some finite trace $\tau_0$
on~$A$, then $\phi$ has a canonical extension to a strictly
continuous state on $B(\F(X))$. We shall see later in
\thmref{old1.3} that, in general, the existence of $\tau_0$ is
also necessary for $\phi$ to have a strictly continuous state
extension to $B(\F(X))$. If $\tilde \phi$ is
the strongly continuous extension of $\phi$ to $B(\F(X))$, then
$\tau_0 (a)=\tilde \phi( P_0 a P_0)$, where $P_0$ is the
projection onto the 'vacuum', i.e. the zeroth component $X^0 = A
\subset \F(X)$.
%since any $(\gamma,\beta)$-KMS state
%on $K(\F(X))$ is of the form $\Ttr_{\tau_0}(\cdot\,\Gamma(e^{-\beta
%D}))$ for some $\tau_0$, by \thmref{old1.3}, and since $\To$ is
%strictly dense in $B(\F(X))$,
\begin{definition}
Let $\phi$ be a $(\gamma, \beta)$-KMS state and set $\tau=
\phi|_A $. Following~\cite{EL2} we say that $\phi$ is of {\em
finite type} if $\tau = \sum^\infty_{n=0}F^n\tau_0$ for some
finite trace $\tau_0$ and we say that $\phi$ is of {\em infinite
type} if $F(\tau)= \tau$. Since $F\tau \leq \tau$ and $F^n\tau(1)
= \Ttr_{\tau}( 1 \, (e^{-\beta D})^{\otimes n})$, we see that
$\phi$ is of finite type iff $\Tr ((e^{-\beta D})^{\otimes n}) \to
0$, and of infinite type iff $\Tr (e^{-\beta D}) =1$.
\end{definition}
Note that in the last part of the proof of \thmref{mainKMSb} we
showed that every $(\gamma, \beta)$-KMS state is a weak* limit of
$(\gamma^\eps, \beta)$-KMS states of finite type as the
perturbation $\eps$ tends to zero.

With the appropriate convention, the above definition makes sense
also for $\beta = \infty$, and it is clear that
KMS$_\infty$-states are necessarily of finite type. As in several
other similar contexts, there is a 'Wold decomposition' for
$(\gamma,\beta)$-KMS states of $\To$.
\begin{proposition}
Under the assumptions of \thmref{old3.1} let $\phi$ be a
$(\gamma,\beta)$-KMS state on $\To$. Then there exists a unique
convex decomposition $\phi=\lambda\phi_1+(1-\lambda)\phi_2$ such
that $\phi_1$ is a $(\gamma,\beta)$-KMS state of finite type and
$\phi_2$ is a $(\gamma,\beta)$-KMS state of infinite type.
\end{proposition}

\begin{proof}
By \thmref{old3.1} we may carry out the decomposition at the level
of traces on $A$; that is, we prove that for any finite trace
$\tau$ on $A$ such that $F\tau\le\tau$ there exists a unique
decomposition $\tau=\tau_1+\tau_2 $ with $\tau_1 =
\sum^\infty_{n=0}F^n\tau_0$ and $F\tau_2=\tau_2$. The uniqueness
is obvious, since $\tau_0$ must equal $\tau-F\tau$. To prove
existence, set $\tau_0=\tau-F\tau$ and $\tau_2=\lim_nF^n\tau$ (the
limit exists because $F^{n+1}\tau\le F^n\tau$). Then
$\sum^m_{n=0}F^n\tau_0=\tau-F^{m+1}\tau\to\tau-\tau_2$ weakly.

The only thing left to check is that $F\tau_2=\tau_2$.
Since $\tau_2\le F^n\tau$, we have $F\tau_2\le F^{n+1}\tau$, and
so $F\tau_2\le\tau_2$. For $a\in A_+$ and $\eps>0$ we can
find a finite subset $I$ of $X$ such that $S_I:=\sum_{\xi\in
I}\theta_{\xi,\xi}\le1$ and $\Tr((1-S_I)ae^{-\beta D})<\eps$. Since
$F^n\tau\le\tau$, we have $\Ttr_{F^n\tau}((1-S)ae^{-\beta
D})<\eps$ for every $n$. Since $F^n\tau$ converges to $\tau_2$
weakly, there exists $n$ such that
$$
\Ttr_{F^n\tau}(Sae^{-\beta D}) <\Ttr_{\tau_2}(Sae^{-\beta D})+\eps,
$$
so
\begin{eqnarray*}
(F\tau_2)(a)&\ge&\Ttr_{\tau_2}(Sae^{-\beta D})
>\Ttr_{F^n\tau}(Sae^{-\beta D})-\eps
>\Ttr_{F^n\tau}(ae^{-\beta D})-2\eps\\
&=&(F^{n+1}\tau)(a)-2\eps\ge\tau_2(a)-2\eps,
\end{eqnarray*}
since $\eps$ was arbitrary, this yields $F\tau_2\ge\tau_2$, and
hence $F\tau_2=\tau_2$.
\end{proof}

We now turn our attention to the KMS states of $\O$. Notice that
if $\phi$ is a $(\gamma,\beta)$-KMS state on $\O$ and we compose
it with the quotient map $\To \to \O$, we obtain a
$(\gamma,\beta)$-KMS state on~$\To$. Thus in order to describe the
KMS-states on $\O$ we only need to describe the KMS-states on
$\To$ that vanish on the kernel of the quotient map.

\begin{theorem} \label{3.2}
Under the assumptions of \thmref{old3.1} suppose $\phi$ is a
$(\gamma,\beta)$-KMS state on $\To$ and let $\tau=\phi|_A$. Then
$\phi$ defines a state on $\O$ if and only if $\Tr(ae^{-\beta
D})=\tau(a)$ for every $a\in I_X$ (where $\Tr(ae^{-\infty D}) =
0$, by convention).
\end{theorem}

\begin{proof}
Suppose first $\beta<\infty$. From the proof of \thmref{old3.1} we
know that $\phi\circ j_X=\Tr(\cdot\,e^{-\beta D})$ on $K(X)$.
Thus, $\Tr(ae^{-\beta D})=\tau(a)$ is equivalent to
$\phi(j_X(a))=\phi(a)$ for $a\in I_X$. This is clearly a necessary
condition for $\phi$ to define a state of $\O$. To see that it is
also sufficient, let $P$ be the projection onto the zeroth
component of the Fock module $\F(X)$; then $b-j_X(b)=PbP=bP$ for
$b\in I_X$ so that $(a-j_X(a))^*(a-j_X(a))=a^*a-j_X(a^*a)$ and
hence $\phi((a-j_X(a))^*(a-j_X(a))) = 0$ for every $a\in I_X$.
Since $\phi$ is a KMS state, the set
$N=\{x\in\To\,|\,\phi(x^*x)=0\}$ is a two-sided ideal in $\To$.
Hence it contains the ideal generated by the elements $a-j_X(a)$
with $a\in I_X$, which is, by definition, the kernel of the map
$\To\to\O$.

For $\beta = \infty$, $\phi$ is a generalized Fock state of $\To$,
so $\phi \circ j_X = 0$. If $\phi$ defines a state on~$\O$, then
$\phi(a) = \phi(a - j_X(a)) = 0$ for $a\in I_X$. Conversely,
suppose $\phi$ vanishes on $I_X$. The generalized Fock state
$\phi$ extends to a state on $B(\F(X))$ with the property
$\phi(x)=\phi(PxP)$. Then for any $x,y\in B(\F(X))$ and $a\in I_X$
we get $\phi(x(a-j_X(a))y)=\phi(PxPaPyP)=0$, since $PxPIPyP\subset
PIP$ for any ideal $I$ in $A$ and since $\phi$ vanishes on $I_X$.
\end{proof}

\begin{remark}
If $X$ is finite over $A$, then $I_X = A$, so a
$(\gamma,\beta)$-KMS state on $\To$ gives one on~$\O$ if and only
if it is of infinite type. In this case, the Wold decomposition of
a KMS state corresponds to the usual essential-singular
decomposition of a state relative to the kernel of the quotient
map $\To \to \O$. However, we point out that in the case of
$\mathcal O_\infty$ and many other interesting situations, cf
\cite{FR,EL2}, the ideal $I_X$ is trivial so $\O$ is actually
equal to $\To$.
\end{remark}

Considering the way in which the quasi-free states of the CAR
algebra are determined by their two-point functions, it is natural
to refer to states given by \eqref{2PtFunction} as {\em
quasi--free states} of~$\To$. Following Evans~\cite{Ev}, we wish
to consider next a slightly more general notion of quasi-free
states; specifically, we wish to allow for a different positive
trace-class operator on each tensor factor. Not surprisingly, the
appropriate formula is easy to guess, and the crux of the matter
is to determine sufficient compatibility conditions on the various
ingredients for it to actually define a state.

\begin{proposition}
Let $\{\tau_n\}^\infty_{n=0}$ be a sequence of traces on $A$ such
that $\tau_0$ is a tracial state and let $\{S_n\}^\infty_{n=1}$ be a
sequence of positive operators in $B_A(X)$ such that
$\Ttr_{\tau_n}(aS_n)\le\tau_{n-1}(a)$ for every $a\in A_+$ and
$n\ge1$. Then there exists a unique gauge-invariant state $\phi$
on $\To$ such that $\phi|_A=\tau_0$ and
$$
\phi(T_{\xi_1}\ldots T_{\xi_n}T^*_{\zeta_n}\ldots T^*_{\zeta_1})
=\tau_n(\<\zeta_1\otimes\ldots\otimes\zeta_n,
S_1\xi_1\otimes\ldots\otimes S_n\xi_n\>)
$$
for every $\xi_1,\ldots,\xi_n,\zeta_1,\ldots,\zeta_n\in X$. If, in
addition, $\Ttr_{\tau_n}(aS_n)=\tau_{n-1}(a)$ for every $a\in I_X$
and $n\ge1$, then $\phi$ defines a state on $\O$.
\end{proposition}

\begin{proof}
We shall use an argument borrowed from the proof of
\cite[Proposition 12.6]{EL2}. Let ${\mathcal T}_0$ be the
subalgebra of gauge-invariant elements of $\To$; since we are
looking for a gauge-invariant state, that is, one that factors
through the conditional expectation $E =
\frac{1}{2\pi}\int^{2\pi}_0\gamma_t (\cdot) dt$ of $\To$ onto
${\mathcal T}_0$, it is enough to define $\phi$ on ${\mathcal
T}_0$. Let $I_n$ be the C$^*$-algebra generated by elements $T_\xi
T^*_\eta$ with $\xi,\eta\in X^{\otimes n}$, and let
$A_n=A+I_1+\ldots+I_n$ (where $A_0=I_0=A$). Then $I_n$ is an ideal
in $A_n$ and ${\mathcal T}_0=\overline{\cup_n A_n}$. The
submodules $X^{\otimes n}$ are ${\mathcal T}_0$-invariant, and we
denote by $\pi_n$ the natural representation of ${\mathcal T}_0$
on $X^{\otimes n}$. Then $\pi_{n-1}$ is faithful on $A_{n-1}$ and
zero on $I_n$, so $A_{n-1}\cap I_n=0$. Consider a positive linear
functional $\psi_n$ on $A_n$ defined by
$$
\psi_n(x)=\Ttr_{\tau_n}(\pi_n(x)(S_1\otimes\ldots\otimes S_n))
$$
for $n\ge1$, with $\psi_0=\tau_0$. Since
$\pi_n(x)=\pi_{n-1}(x)\otimes1$ for every $x\in A_{n-1}$, we have
$\psi_n|_{A_{n-1}}\le\psi_{n-1}$ for $n\ge1$ by
Proposition~\ref{1.2}, and we may define selfadjoint linear
functionals $\phi_n$ on $A_n=A_{n-1}\oplus I_n$ by induction:
$\phi_n=\phi_{n-1}\oplus\psi_n$, with $\phi_0=\tau_0$. Since
$\phi_n|_A=\tau_0$ is a state and the left action of $A$ is
non-degenerate, to prove that $\phi_n$ is a state it is enough to
check that it is positive. We shall prove by induction that
$\phi_n\ge\psi_n$. This is true for $n=0$, since $\psi_0=\phi_0$.
If this is true for $n-1$, then
$\psi_n|_{A_{n-1}}\le\psi_{n-1}\le\phi_{n-1}$. For $x\in A_{n-1}$
and $y\in I_n$ we have
$$
\phi_n((x+y)^*(x+y))=\phi_{n-1}(x^*x)+\psi_n(x^*y+y^*x+y^*y)
\ge\psi_n(x^*x+x^*y+y^*x+y^*y),
$$
so $\phi_n\ge\psi_n$. Thus each $\phi_n$ is a state, and since
the sequence is coherent in the sense that
$\phi_n|_{A_{n-1}}=\phi_{n-1}$, there exists a unique state $\phi$
on ${\mathcal T}_0$ such that
$\phi|_{I_n}=\phi_n|_{I_n}=\psi_n|_{I_n}$.

Suppose now that $\Ttr_{\tau_n}(aS_n)=\tau_{n-1}(a)$ for every
$a\in I_X$; we shall prove that $\phi$ is zero on the two-sided
ideal generated by the elements $a-j_X(a)$, for $a\in I_X$. Let
${\mathcal T}_{n,m}$ be the linear span of elements of the form
$T_\xi T^*_\zeta$, with $\xi\in X^{\otimes n}$ and $\zeta\in
X^{\otimes m}$ (where ${\mathcal T}_{0,0}=A$). We must prove that
$\phi(xay)=\phi(xj_X(a)y)$ for $x\in{\mathcal T}_{n,m}$,
$y\in{\mathcal T}_{n',m'}$. Since $aT_\xi=j_X(a)T_\xi$ for every
$\xi\in X$, we have $xay=xj_X(a)y$ if either $n'>0$, or $m>0$.
Thus we may restrict our attention to the case $n'=m=0$. Because
of gauge-invariance we also have $\phi(xay)=\phi(xj_X(a)y)=0$ if
$n\ne m'$. So it remains to consider only the case when $n=m'$,
$m=n'=0$. Let $x=T_{\xi_0}$ and $y=T^*_{\zeta_0}$, with
$\xi_0,\zeta_0\in X^{\otimes n}$. Denoting
$S_1\otimes\ldots\otimes S_n$ by $\tilde S_n$, for any
$\xi,\zeta\in X$ we have
$$
\phi(T_{\xi_0}j_X(\theta_{\xi,\zeta})T^*_{\zeta_0})=\phi(T_{\xi_0}T_\xi
T^*_\zeta T^*_{\zeta_0})=\tau_{n+1}(\<\zeta_0\otimes\zeta,\tilde
S_n\xi_0\otimes S_{n+1}\xi\>) =\Ttr_{\tau_{n+1}}(\<\zeta_0,\tilde
S_n \xi_0\>S_{n+1}\theta_{\xi,\zeta}),
$$
so for any $a\in I_X$ we get
\begin{eqnarray*}
\phi(T_{\xi_0}j_X(a)T^*_{\zeta_0})
&=&\Ttr_{\tau_{n+1}}(\<\zeta_0,\tilde
S_n \xi_0\>S_{n+1}a)=\tau_n(a\<\zeta_0,\tilde S_n\xi_0\>)
=\tau_n(\<\zeta_0a^*,\tilde S_n\xi_0\>)\\
&=&\phi(T_{\xi_0}T^*_{\zeta_0a^*})=\phi(T_{\xi_0}aT^*_{\zeta_0}).
\end{eqnarray*}
\end{proof}

In order to illustrate the range of application of the general
theory we discuss several situations that have appeared in the
literature and which are unified by the present approach. Other
systems such as those studied in \cite{E1,E2,MWY} can be analyzed
in a similar way.

\begin{example}\cite{OP,BEH,Ev}
Let $X=H$ be a Hilbert space considered as a Hilbert bimodule over
$A=\C$. Then $\To$ is the Toeplitz-Cuntz algebra $\mathcal T_n$
and $\O$ is the Cuntz algebra $\mathcal O_n$ corresponding to $n =
\dim H$. If $\tau$ is the unique tracial state on $A$, then $\Tr$
is the usual trace on $B(H)$. Let $\{U_t=e^{itD}\}_{t\in \mathbb
R}$ be a one-parameter unitary group on $H$, and let $\gamma$ be
the corresponding quasi-free dynamics. Then our results say that a
$(\gamma,\beta)$-KMS state on ${\mathcal T}_H$ exists if and only
if $\Ttr(e^{-\beta D})\le1$, and such a state extends to a normal
state on $B(\F(H))$ if and only if $\Ttr(e^{-\beta D})<1$. If
$\dim H\ge2$ and $\beta\ge0$ (or $\beta\le0$), then the condition
$\Ttr(e^{-\beta D})\le1$ implies that $D$ is positive (resp.
negative) and non-singular, so the $(\gamma,\beta)$-KMS state is
unique. If $H$ is infinite-dimensional, then $I_H=0$ and
${\mathcal O}_H= {\mathcal T}_H$. If $\dim H<\infty$, then
$I_H=\C$, so a $(\gamma,\beta)$-KMS state on ${\mathcal O}_H$
exists if and only if $\Ttr(e^{-\beta D})=1$.
\end{example}

\begin{example}\cite{EL2}\label{exampleel2}
Let $T=(T(j,k))_{j,k\in I}$ be a (possibly infinite) $0-1$ matrix
with no identically zero columns and rows. Consider the rows
$q_j=(T(j,k))_{k\in I}$ as elements of $l^\infty(I)$, and let $A$
be the C$^*$-algebra they generate. Let $X$ be the Hilbert
$A$-bimodule generated as a right Hilbert $A$-module by vectors
$\xi_j$, $j\in I$, such that $\<\xi_j,\xi_k\>=\delta_{jk}q_j$,
$q_j\xi_k=T(j,k)\xi_k$. By~\cite{Sz} the algebra $\O$ is the
Cuntz-Krieger algebra corresponding to the matrix $T$ as in
\cite{EL1}. If $F$ is a finite subset of $I$, then $\sum_{j\in
F}\theta_{\xi_j,\xi_j}$ is the projection onto the right submodule
generated by $\xi_j$, $j\in F$. It follows that if $\tau$ is a
trace on~$A$, then $\Tr(S)=\sum_j\tau(\<\xi_j,S\xi_j\>)$ for any
$S\in B(X)_+$.

Let the dynamics on $X$ be given by $U_t\xi_j=N^{it}_j\xi_j$ for
some choice $N_j >1$ for $j = 1, 2, \ldots$, which ensures that
$U$ satisfies the positivity condition. Then $(\gamma,\beta)$-KMS
states on $\To$ are in one-to-one correspondence with states on
$A$ such that $\sum_jN^{-\beta}_j\tau(\<\xi_j,a\xi_j\>)\le\tau(a)$
for any $a\in A_+$. Any positive element in $A$ can be
approximated by a linear combination with positive coefficients of
projections $q(Y,Z)=\prod_{l\in Y}q_l\prod_{k\in Z}(1-q_k)$, where
$Y$ and $Z$ are finite subsets of $I$. Note that
$q(Y,Z)\xi_j=T(Y,Z,j)\xi_j$, where
$$
T(Y,Z,j)=\prod_{l\in Y}T(l,j)\prod_{k\in Z}(1-T(k,j)).
$$

Thus $(\gamma,\beta)$-KMS states on $\To$ correspond to states
$\tau$ on $A$ such that
\begin{equation} \label{e4.10}
\sum_jN^{-\beta}_jT(Y,Z,j)\tau(q_j)\le\tau(q(Y,Z))\ \ \hbox{for
all finite}\ Y,Z\subset I,
\end{equation}
namely, to $\beta$-subinvariant states of $A$ as defined in
\cite[12.3]{EL2}. Notice that our result bypasses the intermediate
step of having to consider measures on path space that are
rescaled by the partial action and goes straight to the
coefficient algebra. The linear span of projections $q(Y,Z)$ such
that the function $I\ni j\mapsto T(Y,Z,j)$ has finite support is
dense in $I_X$. So to obtain a state on $\O$ the inequality
\eqref{e4.10} must be equality for all $Y$ and $Z$ such that the
function $T(Y,Z,\cdot)$ has finite support. We point out that the
dynamics arising from such infinite matrices do not satisfy in
general the fullness assumption of \cite{PWY}.
\end{example}

\begin{example}\cite{BEH,BEK} \label{examplebeh}
Let $A$ be a unital C$^*$-algebra, $p$ a full projection in $A$,
$\alpha$ an injective $*$-endomorphism on $A$ such that
$\alpha(A)=pAp$. Consider the semigroup crossed product C*-algebra
$A\rtimes_\alpha\NN$. Let~$\gamma$ be the periodic dynamics on
$A\rtimes_\alpha\NN$ defined by the dual action of ${\mathbb T}$.
It is known that $A\rtimes_\alpha\NN$ can be considered as a
Cuntz-Krieger-Pimsner algebra: the module $X$ is the space $Ap$
with left and right actions of $A$ given by $a\cdot\xi\cdot
b=a\xi\alpha(b)$, and with inner product
$\<\xi,\zeta\>=\alpha^{-1}(\xi^*\zeta)$. Then $\gamma$ is the
gauge action. Since $p$ is full, the module is finite, so $I_X=A$.
Thus $(\gamma,\beta)$-KMS states on $\To$ (resp. $\O$) correspond
to traces $\tau$ such that $e^{-\beta}\Tr(a)\le\tau(a)$ (resp.
$e^{-\beta}\Tr(a)=\tau(a)$) for $a\in A_+$. Since $p$ is full, to
prove an equality/inequality for traces on $A$ it is enough to
check it on $pAp=\alpha(A)$. Noting that $p=\theta_{p,p}$ in
$B(X)$, for any $a\in A$ we get
$\Tr(\alpha(a))=\Tr(\alpha(a)\theta_{p,p})=\tau(\<p,\alpha(a)p\>)
=\tau(a)$. Thus $(\gamma,\beta)$-KMS states on $\To$ (resp. $\O$)
correspond to traces $\tau$ on $A$ such that
$e^{-\beta}\tau\le\tau\circ\alpha$ (resp.
$e^{-\beta}\tau=\tau\circ\alpha$). It is proved in \cite{BEH} that
any closed subset of $(0,+\infty)$ can be realized as the set of
possible temperatures of the system $(\O,\gamma)$ for a convenient
choice of AF-algebra $A$ and endomorphism~$\alpha$ with $\O$ is
simple. In \cite{BEK} a similar construction yields quite
arbitrary choices of the simplex of $\beta$-KMS states for each
$\beta$.
\end{example}

\begin{remark}\label{cooling}
A question raised in the introduction of \cite{EL2} is whether the
Toeplitz Cuntz-Krieger algebra of an infinite matrix can have KMS
states of finite and of infinite type coexisting at a given finite
inverse temperature. We would like to shed some light on the
analogous question for the quasi--free dynamics on the
Toeplitz--Pimsner algebras. Under the same assumptions of
\thmref{mainKMSb}, let $\beta_0 < \infty$ and suppose that $\phi$
is a $(\gamma,\beta_0)$-KMS state with restriction $\tau$ to $A$.
For each $\beta > \beta_0$, we then have that $F_\beta \tau \leq
F_{\beta_0} \tau \leq \tau $ so $\tau$ determines a
$(\gamma,\beta)$-KMS state $\phi_\beta$ of $\To$, which cannot be
of infinite type because of our positivity assumption. The
parameters of Example \ref{examplebeh} can be adjusted to get also
a $(\gamma, \beta)$-KMS state of $\O$, and the resulting infinite
type $(\gamma, \beta)$-KMS state of $\To$ will thus coexist,
 at inverse temperature $\beta$, with
the finite part of the decomposition of $\phi_\beta$. A further
strengthening of the positivity assumption, namely, the assumption
that $D \geq c 1$ for some $c >0$, yields $F_\beta^n \tau \leq
e^{-n(\beta - \beta_0)c} F_{\beta_0}^n \tau$, from which one
readily sees that the state $\phi_\beta$ above is of finite type
for $\beta >\beta_0$. This gives a Pimsner algebra version of the
'cooling lemma' \cite[Lemma 15.1]{EL2}, which implies that KMS
states of finite type are weak* dense in all KMS states.
\end{remark}

\section{KMS weights on module algebras} \label{kmsweights}
The theory of KMS weights on C$^*$-algebras is similar to (and is
based on) the theory of normal weights on von Neumann algebras. We
refer the reader to~\cite{S} and~\cite{K} for the basic definitions.
Suppose $\R\ni t\mapsto\sigma_t$ is a one-parameter automorphism
group of a C$^*$-algebra $A$. We assume that $\sigma$ is
continuous in the sense that the function $\R\ni
t\mapsto\sigma_t(a)$ is norm-continuous for all $a\in A$. The same
assumption is made for more general one-parameter group of
isometries on Banach spaces. In several places where we consider
von Neumann algebras the continuity assumption is weaker: the
function $\R\ni t\mapsto\sigma_t(a)$ is weakly (operator) continuous.

Let $\phi$ be a weight on $A$. As usual, we set
\[
\M^+_\phi=\{a\in A_+\,|\,\phi(a)<\infty\},\quad
\N_\phi=\{a\in A\,|\,a^*a\in\M^+_\phi\}, \text{ and} \quad
\M_\phi=\span\,\M^+_\phi=\N^*_\phi\N_\phi,
\]
and extend $\phi$ to a linear functional on $\M_\phi$. We say that
$\phi$ is a $(\sigma,\beta)$-KMS weight if $\phi$ is lower
semicontinuous on $A_+$, densely defined (i.e. $\M^+_\phi$ is
dense in $A_+$), $\sigma$-invariant, and satisfies the
KMS-condition in the form
\[
\phi(x^*x)=\phi(\sigma\sb(x)\sigma\sb(x)^*)\quad \text{ for every }
x\in D(\sigma\sb),
\]
where $D(\sigma\sb)$ is the domain of definition of $\sigma\sb$.

Given such a weight $\phi$, the well known GNS construction
produces a Hilbert space $H_\phi$ and a linear map
$\Lambda_\phi\colon\N_\phi\to H_\phi$ such that
$\Lambda_\phi(\N_\phi)$ is dense in $H_\phi$ and
$(\Lambda_\phi(x),\Lambda_\phi(y))=\phi(y^*x)$. There is a
representation $\pi_\phi$ of $A$ on $H_\phi$, defined by letting
$\pi_\phi(x)\Lambda_\phi(y)=\Lambda_\phi(xy)$, and the set
$\U_\phi=\Lambda_\phi(\N_\phi\cap\N^*_\phi)$ is a left Hilbert
algebra with operations
\[
\Lambda_\phi(x)\Lambda_\phi(y)=\Lambda_\phi(xy),\ \
\Lambda_\phi(x)^\#=\Lambda_\phi(x^*).
\]
The associated von Neumann algebra $\L(\U)=\pi_\phi(A)''$ is
equipped with a canonical normal semifinite faithful (n.s.f.)
weight $\Phi$. Then $\phi=\Phi\circ\pi_\phi$ and
$\sigma^\Phi_t\circ\pi_\phi=\pi_\phi\circ\sigma_{-\beta t}$, where
$\sigma^\Phi_t$ is the modular group of~$\Phi$.

Any lower semicontinuous densely defined weight on $A$ extends
uniquely to a strictly lower semicontinuous weight $\bar\phi$ on
the multiplier algebra $M(A)$. If $\phi$ is a KMS-weight, then
$\bar\phi=\Phi\circ\bar\pi_\phi$, where $\bar\pi_\phi$ is the
canonical extension of $\pi_\phi$ to a representation of $M(A)$.
The weight $\bar\phi$ is only strictly densely defined, and the
one-parameter automorphism group $\sigma$ is only strictly
continuous on $M(A)$, so $\bar\phi$ is not a KMS-weight in the
sense of the definition above. However, since
$\bar\phi=\Phi\circ\bar\pi_\phi$, $\bar\phi$ satisfies a form of
the KMS-condition, so that its restriction to a $\sigma$-invariant
subalgebra $B$ of $M(A)$ is a KMS-weight if $\sigma|_B$ is
continuous and $\bar\phi|_B$ is densely defined.

\smallskip

We now consider a right Hilbert $A$-module $X$ and a continuous
one-parameter group of isometries $\R\ni t\mapsto U_t$ of $X$ such
that $\<U_t\xi,U_t\zeta\>=\sigma_t(\<\xi,\zeta\>)$. Then $U_t(\xi
a)=(U_t\xi)\sigma_t(a)$, and we define a dynamics~$\gamma = \Ad U$
on $K(X)$ by $\gamma_t(T)=U_tTU_{-t}$. We shall extend
\thmref{1.1} to the present situation, for which we need to
associate an induced weight $\kappa_\phi$  on $B(X)$ to each
$(\sigma,\beta)$-KMS weight $\phi$ on $A$. We give two equivalent
constructions of this induced weight.

\smallskip

In the first construction, we replace $A$ by $\overline{\<X,X\>}$
so we may assume that $X$ is full. Let $\pi\colon A\to B(H)$ be an
arbitrary representation of $A$ such that there exists a n.s.f.
weight $\Phi$ on $L=\pi(A)''$ with $\phi=\Phi\circ\pi$ and
$\sigma^\Phi_t\circ\pi=\pi\circ\sigma_{-\beta t}$. Set
$M=\pi(A)'$. Consider the Hilbert space $H_X=X\otimes_A H$ and the
induced representation $\rho$ of $B(X)$ on $H_X$. Set
$N=\rho(B(X))''$. Note that $M$ acts faithfully on~$H_X$, and
$N'=M$ in $B(H_X)$. Let now $\Phi'$ be an arbitrary n.s.f. weight
on $M$ and $\Psi$ the unique n.s.f. weight on $N$ such that
\[
\Delta(\Psi/\Phi')^{it}=U_{-\beta t}\otimes\Delta(\Phi/\Phi')^{it}\ \
\hbox{in}\ \ B(H_X),
\]
where $\Delta(\cdot/\cdot)$ denotes the spatial derivative, cf
\cite{Co2,S}. We then set $\kappa_\phi=\Psi\circ\rho$.

Let us compute $\kappa_\phi$ explicitly on a dense subalgebra of
$K(X)$. Let $\zeta\in H$ be a $\Phi'$-bounded vector, that is, the
map $R_\zeta\colon H_{\Phi'}\to H$, given by
$R_\zeta(\Lambda_{\Phi'}(x))=x\zeta$ for $x\in\N_{\Phi'}$, is
bounded. For any $\xi\in X$ the vector $\xi\otimes\zeta\in H_X$ is
$\Phi'$-bounded as well and $R_{\xi\otimes\zeta}\eta=\xi\otimes
R_\zeta\eta$ for $\eta\in H_{\Phi'}$. Then if $\xi\in D(U\sa)$ and
$\zeta\in D(\Delta(\Phi/\Phi')^\2)$ is $\Phi'$-bounded, we get by
definition
\begin{eqnarray*}
\Psi(R_{\xi\otimes\zeta}R^*_{\xi\otimes\zeta})
&=&\|\Delta(\Psi/\Phi')^\2(\xi\otimes\zeta)\|^2
  =\|U\sa\xi\otimes\Delta(\Phi/\Phi')^\2\zeta\|^2\\
&=&(\Delta(\Phi/\Phi')^\2\zeta,\pi(\<U\sa\xi,U\sa\xi\>)
\Delta(\Phi/\Phi')^\2\zeta).
\end{eqnarray*}
At this point it is convenient to introduce the form
$\Phi(x\sigma^\Phi\sd(y))$, whose relevant properties are stated
in the following essentially known lemma.

\begin{lemma} \label{old1.5}
Let $\Phi$ be a n.s.f. weight on a von Neumann algebra $M$,
with modular group $\sigma$. Then the bilinear form
$(\cdot,\cdot)_\Phi\colon\M_\Phi\times D(\sigma\sd)\to\C$, defined by
$(x,y)_\Phi=\Phi(x\sigma\sd(y))$, extends to a bilinear form on
$\M_\Phi\times M$ with the following properties:
\begin{enumerate}
\item[(i)] for $x\in\M_\Phi$, $(x,\cdot)_\Phi$ is a normal
linear functional on $M$; if $x\in\M^+_\Phi$, then
$(x,\cdot)_\Phi$ is positive and $(x,\cdot)_\Phi\le\|x\|\Phi$;

\item[(ii)] for $x,y\in\M_\Phi$, $(x,y)_\Phi=(y,x)_\Phi$;

\item[(iii)] for $x_1,x_2\in\N_\Phi\cap D(\sigma\sc)$,
$(x^*_2x_1,\cdot)_\Phi=\Phi(\sigma\sc(x_1)\cdot\sigma\sc(x_2)^*)$;

\item[(iv)] if $\{x_k\}_k$ is a net in $\M^+_\Phi\cap
D(\sigma\sd)$ such that $x_k\le1$, the set $\{\sigma\sd(x_k)\}_k$
is bounded in norm, and $\sigma\sd(x_k)\to1$ strongly, then
$\Phi(x)=\lim_k(x_k,x)_\Phi$ for any $x\in M_+$.
\end{enumerate}
\end{lemma}

\begin{proof} Consider the GNS-representation
$\Lambda\colon\N_\Phi\to H$ as described above. Recall that, in
terms of the modular group $\sigma$, the modular involution $J$ is
given by $J\Lambda(x)=\Lambda(\sigma\sc(x)^*)$ for
$x\in\N_\Phi\cap D(\sigma\sc)$. It has the properties
\begin{eqnarray*}
Jy^*J\Lambda(x) &=& \Lambda(x\sigma\sd(y))\quad
\forall x\in\N_\Phi \quad \forall y\in D(\sigma\sd), \text{ and}\\
 JyJ\Lambda(x) &=& xJ \Lambda(y) \quad \forall x,y\in\N_\Phi.
\end{eqnarray*}
Now if $x=x^*_2x_1$, $x_1,x_2\in\N_\Phi$, $y\in D(\sigma\sd)$,
then
\[
(x,y)_\Phi=(x^*_2x_1,y)_\Phi=(\Lambda(x_1\sigma\sd(y)),\Lambda(x_2))
=(Jy^*J\Lambda(x_1),\Lambda(x_2)).
\]
This shows that $(x,\cdot)_\Phi$ extends to a normal linear
functional on $M$. Moreover, if $x\ge0$, then we can take
$x_1=x_2=x^\2$ and conclude that $(x,y)_\Phi\ge0$ for $y\ge0$. If
$x,y\in\M^+_\Phi$,
\begin{eqnarray*}
(x,y)_\Phi&=&(JyJ\Lambda(x^\2),\Lambda(x^\2))
=(Jy^\2J\Lambda(x^\2),Jy^\2J\Lambda(x^\2))
=(x^\2J\Lambda(y^\2),x^\2J\Lambda(y^\2))\\
&=&(JxJ\Lambda(y^\2),\Lambda(y^\2))
\le\|x\|(\Lambda(y^\2),\Lambda(y^\2))=\|x\|\Phi(y).
\end{eqnarray*}
Thus part (i) is proved. Part (ii) is already proved for positive
$x$ and $y$, the general case follows by linearity. If
$x_2,x_1\in\N_\phi\cap D(\sigma\sc)$, we have
\begin{eqnarray*}
(x^*_2x_1,y)_\Phi&=&(Jy^*J\Lambda(x_1),\Lambda(x_2))
=(yJ\Lambda(x_2),J\Lambda(x_1))\\
&=&(y\Lambda(\sigma\sc(x_2)^*),\Lambda(\sigma\sc(x_1)^*))
=\Phi(\sigma\sc(x_1)y\sigma\sc(x_2)^*),
\end{eqnarray*}
which proves (iii). It remains to prove (iv). Let $x\in M_+$.
Since $\{\sigma\sd(x_k)^*x\sigma\sd(x_k)\}_k$ converges weakly to
$x$, and $x^2_k\le x_k$, we have
$$
\Phi(x)\le\liminf_k\Phi(\sigma\sd(x_k)^*x\sigma\sd(x_k))
=\liminf_k(x^2_k,x)_\Phi\le\liminf_k(x_k,x)_\Phi.
$$
On the other hand, $(x_k,x)_\Phi\le\Phi(x)$ by~(i), hence
$\Phi(x)=\lim_k(x_k,x)_\Phi$. This finishes the proof of \lemref{old1.5}.
\end{proof}

Returning to the computation of $\Psi$, for any $a\in L$ we get
$$
(\Delta(\Phi/\Phi')^\2\zeta,a\Delta(\Phi/\Phi')^\2\zeta)
=(\Delta(\Phi/\Phi')^\2\zeta,\Delta(\Phi/\Phi')^\2\sigma^\Phi\sc(a)\zeta)
=\Phi(R_\zeta R^*_{\sigma^\Phi\sc(a)\zeta})
=(R_\zeta R^*_\zeta,a^*)_\Phi.
$$
Thus
\begin{equation} \label{e3.01}
\Psi(R_{\xi\otimes\zeta}R^*_{\xi\otimes\zeta})
=(R_\zeta R^*_\zeta,\pi(\<U\sa\xi,U\sa\xi\>))_\Phi.
\end{equation}
Note that
$R^*_{\xi\otimes\zeta}(\xi_0\otimes\zeta_0)
=R^*_\zeta\pi(\<\xi,\xi_0\>)\zeta_0$.
So if we introduce an operator $T_{\xi,a}$ on $H_X$ defined by
$T_{\xi,a}(\xi_0\otimes\zeta_0)=\xi\otimes
a\pi(\<\xi,\xi_0\>)\zeta_0$, then
$R_{\xi\otimes\zeta}R^*_{\xi\otimes\zeta}=T_{\xi,R_\zeta
R^*_\zeta}$. Then \eqref{e3.01} implies that
\begin{equation*} \label{e3.02}
\Psi(T_{\xi,a})=(a,\pi(\<U\sa\xi,U\sa\xi\>))_\Phi
\end{equation*}
for any $a$ in the algebra $\span\{R_{\zeta_1}R^*_{\zeta_2}\,|\,
\zeta_i\in D(\Delta(\Phi/\Phi')^\2)\ \hbox{is}\
\Phi'\hbox{-bounded},\ i=1,2\}$. If we now apply this to an
approximate unit $a_i\nearrow1$ and use Lemma~\ref{old1.5}(i-ii)
together with the property
$T_{\xi,a_i}\nearrow\rho(\theta_{\xi,\xi})$, we conclude that
\begin{equation} \label{e3.03}
\kappa_\phi(\theta_{\xi,\xi})=\phi(\<U\sa\xi,U\sa\xi\>)
\end{equation}
for every $\xi\in D(U\sa)$ such that $\<U\sa\xi,U\sa\xi\>\in\M_\phi$.

\smallskip

In our second construction of $\kappa_\phi$, we shall show
directly that there is a linear functional satisfying property
\eqref{e3.03} on a dense subalgebra of $K(X)$, and then extend it
to a weight on the whole algebra using the GNS-representation. For
this we introduce the following sets:
\begin{eqnarray*}
C_\phi&=&\{a\in A\,|\,a\ \hbox{is}\ \sigma\text{-analytic and }
\sigma_z(a)\in\N_\phi\cap\N^*_\phi\
\forall z\in\C\},\\
X_0&=&\{\xi\in X\,|\,\xi\ \hbox{is}\ U\text{-analytic}\},\\
X_\phi&=&X_0C_\phi,\\
\U&=&\span\{\theta_{\xi,\zeta}\,|\,\xi,\zeta\in X_\phi\}.
\end{eqnarray*}
Note that if $\xi\in D(U_z)$ and $\zeta\in D(U_{\bar z})$, then
$\theta_{\xi,\zeta}\in D(\gamma_z)$ and
$\gamma_z(\theta_{\xi,\zeta})=\theta_{U_z\xi,U_{\bar z}\zeta}$;
moreover, at the level of the coefficient algebra,
$\<\zeta,\xi\>\in D(\sigma_z)$ and
$\sigma_z(\<\zeta,\xi\>)=\<U_{\bar z}\zeta,U_z\xi\>$. Thus
$C_\phi$ is a dense $*$-subalgebra of $A$, $X_\phi$ is a dense
subspace of $X$, and $\U$ a dense $*$-subalgebra of $K(X)$
consisting of $\gamma$-analytic elements. Choose an approximate
unit $\{e_k=\sum_{\eta\in I_k}\theta_{\eta,\eta}\}_k$ in $\U$ and
define
\begin{equation} \label{e1.1}
\kappa(x)=\lim_k\sum_{\eta\in I_k}\phi(\<U\sa\eta,xU\sa\eta\>).
\end{equation}
for $x\in\U$. If $\eta,\xi\in D(U\sa)$ then
\begin{equation} \label{e1.2}
\phi(\<U\sa\eta,\theta_{\xi,\xi}U\sa\eta\>)
=\phi(\<U\sa\eta,\xi\>\<\xi,U\sa\eta\>)
=\phi(\<U\sa\xi,\eta\>\<\eta,U\sa\xi\>)
=\phi(\<U\sa\xi,\theta_{\eta,\eta}U\sa\xi\>),
\end{equation}
since $\<\xi,U\sa\eta\>\in D(\sigma\sb)$ and
$\sigma\sb(\<\xi,U\sa\eta\>)=\<U\sa\xi,\eta\>$. Hence for $\xi\in
X_\phi$
\[
\kappa(\theta_{\xi,\xi})=\lim_k\sum_{\eta\in
I_k}\phi(\<U\sa\eta,\theta_{\xi,\xi}U\sa\eta\>)
=\lim_k\phi(\<U\sa\xi,e_kU\sa\xi\>)=\phi(\<U\sa\xi,U\sa\xi\>),
\]
so $\kappa$ satisfies \eqref{e3.03}.
It is easy to check that $\kappa$ has the following properties:
\begin{list}{}{}

\item{(i)} $\kappa$ is $\gamma_z$-invariant for any $z\in\C$;

\item{(ii)} $\kappa(xy)=\kappa(\gamma\sb(y)\gamma\sa(x))$ for any
$x,y\in\U$;

\item{(iii)} $\kappa(\theta_{\xi
a,\xi})=(\<U\sa\xi,U\sa\xi\>,a)_\phi$ for any $\xi\in X_\phi$ and
$a\in C_\phi$.

\item{(iv)} the function $z\mapsto\kappa(\gamma_z(x)y)$ is
analytic for any $x,y\in\U$.
\end{list}
Since $\kappa(x^*x)\ge0$ by \eqref{e1.1}, there exist a Hilbert
space $H$ and a linear map $\Lambda\colon\U\to H$ such that
$\Lambda(\U)$ is dense in $H$, and
$(\Lambda(x),\Lambda(y))=\kappa(y^*x)$. The kernel of this map is
the set
$$
\Ker\,\Lambda = \{x\in\U\,|\,\kappa(x^*x)=0\}=\{x\in\U\,|\,
\kappa(yx)=0\ \forall y\in\U\},
$$
which is obviously a $\gamma_z$-invariant left ideal in $\U$.
Property (ii) of $\kappa$ implies that if $x$ is in this ideal,
then also $\gamma\sa(x^*)$, and hence also $x^*$, is in this
ideal. Thus the ideal is self-adjoint and two-sided, and it follows
that $\Lambda(\U)$ has the canonical structure of an algebra with
involution, in particular, $\Lambda(x)\Lambda(y)=\Lambda(xy)$, and
$\Lambda(x)^\#=\Lambda(x^*)$.

Let us check that $\Lambda(\U)$ is a left Hilbert algebra. It is
obvious that $(\Lambda(x)\Lambda(y),\Lambda(z))
=(\Lambda(y),\Lambda(x)^\#\Lambda(z))$. The map
$\Lambda(y)\mapsto\Lambda(x)\Lambda(y)$ is a bounded map for each
$x\in\U$ by (\ref{e1.1}). In fact, we see that the norm of the
mapping is not bigger than~$\|x\|$. Since $\kappa$ is
$\gamma_z$-invariant, property (ii) of $\kappa$ can be rewritten
as
\begin{equation} \label{e1.5}
(\Lambda(x^*),\Lambda(y))
=(\Lambda(\gamma_{i\beta}(y^*)),\Lambda(x)).
\end{equation}
It follows $\Lambda(x)\mapsto\Lambda(x^*)$ is a closable operator.
It remains to prove that $\Lambda(\U)^2$ is dense in
$\Lambda(\U)$. Let $\{e_k\}_k$ be an approximate unit in $\U$. We
claim that $\Lambda(e_k)\Lambda(x)\to\Lambda(x)$ for any $x\in\U$.
Let $x=\theta_{\xi,\zeta}$. Then by property (iii) of $\kappa$
\[
||\Lambda(x)-\Lambda(e_k)\Lambda(x)||^2
=||\Lambda(\theta_{\xi-e_k\xi,\zeta})||^2
=\kappa(\theta_{\zeta\<\xi-e_k\xi,\xi-e_k\xi\>,\zeta})
=(\<U\sa\zeta,U\sa\zeta\>,\<\xi-e_k\xi,\xi-e_k\xi\>)_\phi.
\]
By \lemref{old1.5}, $(\<U\sa\zeta,U\sa\zeta\>,\cdot)_\phi$ is a
bounded linear functional on $A$, so $\{\Lambda(e_kx)\}_k$
converges to $\Lambda(x)$. This proves the claim and completes the
proof that $\Lambda(\U)$ is a left Hilbert algebra.

Let $\Psi$ be the canonical n.s.f. weight on the associated von
Neumann algebra $\L(\U)$. For each $x\in\U$ let $\rho(x)$ be the
operator of left multiplication by $\Lambda(x)$. We have already
shown that $||\rho(x)||\le||x||$. Hence $\rho$ extends by continuity,
first to a representation of $K(X)$, and then to a representation of
$B(X)$. Finally, we define $\kappa_\phi=\Psi\circ\rho$.

Let us compute the modular group of $\Psi$. By property (iv) of
$\kappa$ the vector-valued function $\C\ni
z\mapsto\Lambda(\gamma_z(x))$ is analytic. Hence there exists a
non-singular positive operator $\Delta$ on $H$ such that
$\Delta^z\Lambda(x)=\Lambda(\gamma_{i\beta z}(x))$. Let $J$ be the
anti-linear involution defined by
$J\Lambda(x)=\Lambda(\gamma\sb(x))^\#$. Then
$J\Delta^\2\Lambda(x)=\Lambda(x)^\#$, from which we conclude that
$\Delta_\Psi=\Delta$. Thus
$\sigma^\Psi_t\circ\rho=\rho\circ\gamma_{-\beta t}$. Since
$\kappa_\phi(y^*x)=(\Lambda(x),\Lambda(y))=\kappa(y^*x)$ by
definition of $\Psi$, we see that $\kappa_\phi|_{K(X)}$ is a
$(\gamma,\beta)$-KMS weight such that property (\ref{e3.03}) is
satisfied for all $\xi\in X_\phi\<X_\phi,X_\phi\>$.

\smallskip

We can now state and prove the following generalization of
\thmref{1.1} for KMS weights. In particular, we will show that
both constructions give the same weight.

\begin{theorem} \label{old1.3}
Let $\phi$ be a $(\sigma,\beta)$-KMS weight on $A$. For $T\in
B(X)$, $T\ge0$, set
\begin{equation}\label{kappadef}
\kappa_\phi(T)=\sup\sum_{\xi\in I}\phi(\<U\sa\xi,TU\sa\xi\>),
\end{equation}
where the supremum is taken over all finite subsets $I$ of
$D(U\sa)$ such that $\<U\sa\xi,U\sa\xi\>\in\M^+_\phi$ for every
$\xi\in I$ and $\sum_{\xi\in I}\theta_{\xi,\xi}\le1$. Then
\begin{enumerate}
\item[(i)] $\kappa_\phi|_{K(X)}$ is a $(\gamma,\beta)$-KMS
weight on $K(X)$, and $\kappa_\phi$ is its strictly lower
semicontinuous extension to $B(X)$;

\smallskip
\item[(ii)] there exists an approximate unit $\{e_k=\sum_{\xi\in
I_k}\theta_{\xi,\xi}\}_k$ in $K(X)$ such that $\xi\in D(U\sa)\cap
D(U\sb)$ and $\<U\sa\xi,U\sa\xi\>\in\M^+_\phi$ for every $\xi\in
I_k$, the net $\{\gamma\sa(e_k)\}_k$ is bounded in norm and
converges strictly to $1$; for any such approximate unit we have
\[
\kappa_\phi(T)=\lim_k\sum_{\xi\in I_k}\phi(\<U\sa\xi,TU\sa\xi\>) \ \
\ \ \hbox{for each}\ T\ge0;
\]

\smallskip
\item[(iii)] if $\xi\in D(U\sa)$, then
$\kappa_\phi(\theta_{\xi,\xi})=\phi(\<U\sa\xi,U\sa\xi\>)$, and if this
number is finite, then
\[
(\theta_{\xi,\xi},\cdot)_{\kappa_{\phi}}
=\phi(\<U\sa\xi,\cdot\,U\sa\xi\>);
\]
moreover, if $\tilde X_\phi\subset X_\phi$ is a dense
$U_z$-invariant $C_\phi$-submodule, then $\kappa_\phi|_{K(X)}$ is
the unique $(\gamma,\beta)$-KMS weight such that
$\kappa_\phi(\theta_{\xi,\xi})=\phi(\<U\sa\xi,U\sa\xi\>)$ for
$\xi\in\tilde X_\phi;$

\smallskip
\item[(iv)] if the module is full, then the mapping $\phi\mapsto
\kappa_\phi|_{K(X)}$ defines a one-to-one correspondence between
$(\sigma,\beta)$-KMS weights on $A$ and $(\gamma,\beta)$-KMS
weights on $K(X)$.
\end{enumerate}
\end{theorem}

\begin{proof}
First note that an approximate unit with the properties stated in
part (ii) always exists. Moreover, if $\tilde X_\phi\subset
X_\phi$ is a dense $U_z$-invariant $C_\phi$-submodule, then such
an approximate unit can be found in the algebra
$\tilde\U=\span\{\theta_{\xi,\zeta}\,|\,\xi,\zeta\in\tilde
X_\phi\}$. Indeed, let $\{f_k=\sum_{\eta\in
I_k}\theta_{\eta,\eta}\}_k$ be an approximate unit in $\tilde\U$.
It is well-known that if we set $\tilde e_k=\frac{1}{\sqrt{\pi}}
\int_\R e^{-t^2}\gamma_t(f_k)dt$, then $\{\tilde e_k\}_k$ is an
approximate unit, the net $\{\gamma_z(e_k)\}_k$ is bounded and
converges strictly to~$1$ for any $z\in\C$. But since each $f_k$
is already $\gamma$-analytic, we can replace the integral by a
finite sum such that the element $e_k$, which we thus obtain, is
arbitrarily close (in norm) to $\tilde e_k$, while
$\gamma\sa(e_k)$ is close to $\gamma\sa(\tilde e_k)$.

We have already shown that there exists a strictly lower
semicontinuous weight $\kappa$ on $B(X)$ such that
$\kappa|_{K(X)}$ is $(\gamma,\beta)$-KMS and
$\kappa(\theta_{\xi,\xi})=\phi(\<U\sa\xi,U\sa\xi\>)$ for
$\xi\in\tilde X_\phi$. To prove the theorem it is enough to show that
$\kappa$ satisfies (iii). Indeed, by \lemref{old1.5}(i),(iv),
$\kappa$ can then be defined as in part (ii) and
\eqref{kappadef}. In particular,
$\kappa_\phi=\kappa$, so $\kappa_\phi$ has properties (i)-(iii).
Part (iv) follows then from the uniqueness result in (iii), since
by symmetry for any $(\gamma,\beta)$-KMS weight $\psi$ on $K(X)$
we can define a strictly lower semicontinuous weight $\kappa_\psi$
on $M(A)$ such that $\kappa_\psi|_A$ is a $(\sigma,\beta)$-KMS
weight and
\[
\kappa_\psi(\<\xi,\xi\>)=\psi(\theta_{U\sb\xi,U\sb\xi})\ \ \hbox{for
any}\ \xi\in D(U\sb).
\]

To show (iii), note that for any $\xi,\zeta\in\tilde X_\phi$
\[
(\theta_{\zeta,\zeta},\theta_{\xi,\xi})_\kappa
=\phi(\<U\sa\zeta,\theta_{\xi,\xi}U\sa\zeta\>).
\]
Since both sides above are continuous functions of $\xi$, the
equality holds for any $\zeta\in\tilde X_\phi$ and $\xi\in X$, and,
using \eqref{e1.2}, we obtain that
\begin{equation} \label{e3.04}
(x,\theta_{\xi,\xi})_\kappa=\phi(\<U\sa\xi,xU\sa\xi\>)
\end{equation}
for every $x\in\tilde\U$ and $\xi\in D(U\sa)$. Choosing an
approximate unit in $\tilde\U$ satisfying the conditions of
\lemref{old1.5}(iv) we conclude that
\[
\kappa(\theta_{\xi,\xi})=\phi(\<U\sa\xi,U\sa\xi\>))
\]
for any $\xi\in D(U\sa)$. If this number is finite,
by \lemref{old1.5}(ii) we can rewrite \eqref{e3.04} as
\[
(\theta_{\xi,\xi},x)_\kappa=\phi(\<U\sa\xi,xU\sa\xi\>)
\]
for any $x\in\tilde\U$. Since both sides above are strictly
continuous linear functionals on $B(X)$, the equality holds for
all $x\in B(X)$. By \lemref{old1.5}(iv) the weight $\kappa$
is completely determined by linear functionals
$(\theta_{\xi,\xi},\cdot)_\kappa$ for $\xi\in\tilde X_\phi$, from which the
uniqueness result follows.
\end{proof}

\begin{remark}
\rm\mbox{\ }

(i) In the particular case when $A$ is a full corner $pBp$, $X=Bp$,
$\gamma$ a one-parameter automorphism group of $B$ leaving the
projection $p\in M(B)$ invariant, $\sigma_t=\gamma_t|_A$,
$U_tx=\gamma_t(x)$, the weight $\kappa_\phi|_{K(X)}$ is an
extension of the $(\sigma,\beta)$-KMS weight $\phi$ to a
$(\gamma,\beta)$-KMS weight on $B$, and the map $\psi\mapsto
\kappa_\psi|_A$ going from KMS-weights on $B$ to KMS-weights on
$A$ is just the restriction map. Thus \thmref{old1.3} says that
any $(\sigma,\beta)$-KMS weight on $A$ can be uniquely extended to
a $(\gamma,\beta)$-KMS weight on $B$. Using the linking algebra
the general case could be reduced to this situation. Namely,
$\kappa_\phi|_{K(X)}$ is a unique weight $\psi$ on $K(X)$ such
that $\Phi=\begin{pmatrix}\psi & 0\cr 0 & \phi\end{pmatrix}$ is an
$(\alpha,\beta)$-KMS weight on the linking algebra
$C=\begin{pmatrix}K(X) & X\cr \overline{X} & A\end{pmatrix}$,
where $\alpha_t(\begin{pmatrix}x & \xi\cr \overline{\zeta} &
a\end{pmatrix}) =\begin{pmatrix}\gamma_t(x) & U_t\xi\cr
\overline{U_t\zeta} & \sigma_t(a)\end{pmatrix}$.

(ii) The induction (or extension) results have natural
counterparts for von Neumann algebras, which can be proved by the
same methods or deduced from our results for C$^*$-algebras (see
also~\cite{CZ}).

(iii) The fact that the weight given by our first construction of
the induced weight $\kappa_\phi$ is independent of the choice of
representation~$\pi$ is essentially equivalent to the main result
of~\cite{Sa}.

(iv) The main point in our first construction of the induced
weight $\kappa_\phi$ is an implicit application of Connes' theorem
on existence and uniqueness of a weight with given Radon-Nikodym
cocycle \cite[Theorem 1.2.4]{Co}. On the other hand, our second
construction uses nothing beyond the basic definitions of the
modular theory, and the induction results we have obtained can be
used  in turn to give an alternative proof of Connes' result.
Indeed, let $M$ be a von Neumann algebra, $\phi$ a n.s.f. weight
on $M$ with modular group~$\sigma$, $\R\ni t\mapsto u_t$ a
strongly continuous unitary $1$-cocycle for $\sigma$. Then by our
results (and remark (ii) above) there exists a unique n.s.f.
weight $\Phi$ on $\Mat_2(M)$ with modular group
\[
\sigma^\Phi_t \big(\begin{pmatrix}a & b\cr c &
d\end{pmatrix}\big)=\begin{pmatrix}u_t\sigma_t(a)u^*_t & u_t\sigma_t(b)\cr
\sigma_t(c)u^*_t & \sigma_t(d)\end{pmatrix}
\]
such that $\Phi(\begin{pmatrix}0&0\cr0&d\end{pmatrix})=\phi(d)$
for any $d\in M_+$. Since $p=e_{22}$ is in the centralizer of
$\Phi$, $\Phi(x)=\Phi(pxp)+\Phi((1-p)x(1-p))$ for any $x\in
\Mat_2(M)_+$. So if we set
$\psi(a)=\Phi(\begin{pmatrix}a&0\cr0&0\end{pmatrix})$, then
$\Phi=\begin{pmatrix}\psi&0\cr0&\phi\end{pmatrix}$. Thus $\psi$ is
a n.s.f. weight with Radon-Nikodym cocycle $(D\psi:D\phi)_t=u_t$.
\end{remark}

The following result is an analogue of \proref{1.2} on
induction in stages for KMS-weights. We use the notation
$\kappa^U_\phi$ instead of $\kappa_\phi$ to indicate explicitly
the dynamics used to induce the weight.

\begin{proposition}
Let $\sigma$ (resp. $\gamma$) be a one-parameter automorphism
group of a C$^*$-algebra $A$ (resp. $B$), $X$ a right Hilbert
$A$-module, $Y$ a Hilbert $A$-$B$-module, $U$ (resp. $V$) a
one-parameter group of isometries of $X$ (resp. $Y$) such that
$\<U_t\xi,U_t\zeta\>=\sigma_t(\<\xi,\zeta\>)$,
$\<V_t\xi,V_t\zeta\>=\gamma_t(\<\xi,\zeta\>)$,
$V_ta\xi=\sigma_t(a)V_t\xi$. Let $\phi$ be a $(\gamma,\beta)$-KMS
weight on $B$. Set $\psi=\kappa^V_\phi|_A$. Suppose $\psi$ is
densely defined, so it is a $(\sigma,\beta)$-KMS weight on $A$.
Then
$$
\kappa^{U\otimes V}_\phi(S\otimes1)=\kappa^U_\psi(S)
$$
for any $S\in B(X)$, $S\ge0$.
\end{proposition}

\begin{proof}
We will give a proof based on the properties of spatial
derivatives, but we point out that a proof along the lines of that
of \proref{1.2} is also possible.

Replacing $A$ by $\overline{\<X,X\>}$ and $B$ by
$\overline{\<Y,Y\>}$ we may assume that the modules are full. Let
$\pi\colon B\to B(H)$ be an arbitrary representation of $B$ such
that there exists a n.s.f. weight $\Phi$ on $L=\pi(B)''$ such that
$\phi=\Phi\circ\pi$ and
$\sigma^\Phi_t\circ\pi=\pi\circ\gamma_{-\beta t}$. We then
consider the induced representations of $B(Y)$ on $H_Y=Y\otimes_B H$
and of $B(X\otimes_A Y)$ on $H_{X\otimes Y}=X\otimes_A Y\otimes_B
H$, and we consider the following von Neumann algebras
\begin{enumerate}
\item $M=L'$  in $B(H)$;
\item $N_Y = B(Y)''$, $N_0=A''$, and
$M_0=A'$ in $B(H_Y)$;
\item $N_{X\otimes Y} = B(X\otimes Y)''$,
$N_X=B(X)''$ in $B(H_{X\otimes Y})$.
\end{enumerate}
The algebra $M$ acts faithfully on $H_Y$ and $H_{X\otimes Y}$,
$M_0$ acts faithfully on $H_{X\otimes Y}=(H_Y)_X$, and
\begin{eqnarray*}
N'_Y&=&M\ \ \hbox{in}\ \ B(H_Y);\\
N'_{X\otimes Y}&=&M,\ N'_X=M_0\ \ \hbox{in}\ \ B(H_{X\otimes Y}).
\end{eqnarray*}
Choose a n.s.f. weight $\Phi'$ on $M$. Let $\Psi_Y$ be the n.s.f.
weight on $N_Y$ such that
$$
\Delta(\Psi_Y/\Phi')^{it}=V_{-\beta
t}\otimes\Delta(\Phi/\Phi')^{it}\ \ \hbox{on}\ \ H_Y.
$$
Set $\Psi_0=\Psi_Y|_{N_0}$, so that $\psi=\Psi_0|_A$. Let
$\Phi'_0$, $\Psi_{X\otimes Y}$ and $\Psi_X$ be the n.s.f. weights
on $M_0$, $N_{X\otimes Y}$ and $N_X$, respectively, such that
\begin{eqnarray*}
\Delta(\Psi_0/\Phi'_0)&=&\Delta(\Psi_Y/\Phi')\ \ \hbox{on}\
\ H_Y;\\
\Delta(\Psi_{X\otimes Y}/\Phi')^{it}&=&U_{-\beta t}\otimes
V_{-\beta
t}\otimes\Delta(\Phi/\Phi')^{it}\ \ \hbox{on}\ \ H_{X\otimes Y};\\
\Delta(\Psi_X/\Phi'_0)^{it}&=&U_{-\beta
t}\otimes\Delta(\Psi_0/\Phi'_0)^{it}=\Delta(\Psi_{X\otimes
Y}/\Phi')^{it}\ \ \hbox{on}\ \ H_{X\otimes Y},
\end{eqnarray*}
so that $\Psi_{X\otimes Y}|_{B(X\otimes Y)}=\kappa^{U\otimes
V}_\phi$ and $\Psi_X|_{B(X)}=\kappa^U_\psi$. We have to prove that
$\Psi_{X\otimes Y}|_{N_X}=\Psi_X$. Let $E\colon N_Y\to N_0$ be the
$\Psi_Y$-preserving conditional expectation. By considering $N_Y$
and $N_0$ as subalgebras of $B(H_Y)$ we get an inverse
operator-valued weight $E^{-1}\colon M_0\to M$, cf \cite{H} and
\cite[Corollary 12.11]{S}.
Then by considering $M_0$ and $M$ as subalgebras of $B(H_{X\otimes
Y})$ we get an operator-valued weight $F=(E^{-1})^{-1}\colon
N_{X\otimes Y}\to N_X$. By definition
\[
\Delta(\Psi_0/\Phi'_0)=\Delta(\Psi_Y/\Phi')=\Delta(\Psi_0\circ
E/\Phi')=\Delta(\Psi_0/\Phi'\circ E^{-1})\ \ \hbox{on}\ \ H_Y,
\]
whence $\Phi'\circ E^{-1}=\Phi'_0$. Then
\[
\Delta(\Psi_X\circ F/\Phi')=\Delta(\Psi_X/\Phi'\circ
E^{-1})=\Delta(\Psi_X/\Phi'_0)=\Delta(\Psi_{X\otimes Y}/\Phi')\ \
\hbox{on}\ \ H_{X\otimes Y},
\]
so $\Psi_X\circ F=\Psi_{X\otimes Y}$. But the property of an
operator-valued weight to have a conditional expectation as the
inverse does not depend on the spatial realization,
cf \cite[Theorem 2.2]{Ko}.
Since $E$ is a conditional expectation, we conclude that $F$ is
also a conditional expectation. Hence $\Psi_{X\otimes
Y}|_{N_X}=\Psi_X$.
\end{proof}

The characterization of KMS states of general quasi-free dynamics
is similar to the case of trivial dynamics on the coefficient
algebra, but requires the correspondence just established between
the KMS weights on $A$ and those on $K(X)$.

\begin{theorem}
Let $\sigma$ be a one-parameter automorphism group of a
C$^*$-algebra $A$, $U$ a one-parameter group of isometries of a
Hilbert $A$-bimodule $X$ such that
$\<U_t\xi,U_t\zeta\>=\sigma_t(\<\xi,\zeta\>)$ and
$U_ta\xi=\sigma_t(a)U_t\xi$, and denote by $\gamma$ the
corresponding quasi-free dynamics on the Toeplitz algebra $\To$.
For $\beta\in\R$, let $F$ be the operator mapping
$(\sigma,\beta)$-KMS states of $A$ into weights on $A$, defined by
\[
F\phi=\kappa^U_\phi|_A,
\]
so that $F\phi$ is a $(\sigma,\beta)$-KMS weight on $A$ when it is
densely defined. Then
\begin{enumerate}
\smallskip
\item[(i)] if $\Phi$ is a $(\gamma,\beta)$-KMS state on
$\To$, then $\phi=\Phi|_A$ is a $(\sigma,\beta)$-KMS state on $A$
such that $F\phi\le\phi$;

\smallskip
\item[(ii)] if $\phi$ is a $(\sigma,\beta)$-KMS state on $A$ such
that $F\phi\le\phi$, then there exists a unique gauge-invariant
$(\gamma,\beta)$-KMS state $\Phi$ on $\To$ such that
$\Phi|_A=\phi$; if $\phi=\sum^\infty_{n=0}F^n\phi_0$, then
$\Phi=\nolinebreak\kappa^{\Gamma(U)}_{\phi_0}|_\To$;

\smallskip
\item[(iii)] if $U$ satisfies the 'positive energy' condition
(i.e. the vectors $\xi$ such that $\Sp_U(\xi)\subset(0,+\infty)$
span a dense subspace of $X$), then any $(\gamma,\beta)$-KMS state
of $\To$ is gauge-invariant, so the mapping $\Phi\mapsto\Phi|_A$
defines a one-to-one correspondence between $(\gamma,\beta)$-KMS
states on $\To$ and $(\sigma,\beta)$-KMS states $\phi$ on $A$ such
that $F\phi\le\phi$;
\smallskip

\item[(iv)] a $(\gamma,\beta)$-KMS state $\Phi$ on $\To$
defines a state on $\O$ if and only if $F\phi=\phi$ on $I_X$,
where $\phi=\Phi|_A$.

\end{enumerate}
\end{theorem}

\end{document}